%% file: lamg.tex
\documentclass[final]{siamltex}

\include{styles}
\usepackage{array}

\title{LEAN ALGEBRAIC MULTIGRID (LAMG): \\ FAST GRAPH LAPLACIAN LINEAR SOLVER}
\author{
Oren E. Livne
\thanks{The University of Chicago, Department of Human Genetics, 920 E. 58th St. CLSC 431F, Chicago, IL 60637. Tel: +1-773-702-5898. Email: {\tt livne@uchicago.edu}}
\and
Achi Brandt \thanks{The Weizmann Institute of Science, Department of Mathematics and Computer Science, POB 26 Rehovot 76100, Israel. Tel. +972-8-934-3545. Email: {\tt abrandt@math.ucla.edu}}
}

\begin{document}
\maketitle

\centerline{\it Dedicated to J. Brahms' Symphony No.~1 in C minor, Op.~68}

\begin{abstract}
Laplacian matrices of graphs arise in large-scale computational applications such as semi-supervised machine learning; spectral clustering of images, genetic data and web pages; transportation network flows; electrical resistor circuits; and elliptic partial differential equations discretized on unstructured grids with finite elements. A Lean Algebraic Multigrid (LAMG) solver of the symmetric linear system $Ax=b$ is presented, where $A$ is a graph Laplacian. LAMG's run time and storage are empirically demonstrated to scale linearly with the number of edges.

LAMG consists of a setup phase during which a sequence of increasingly-coarser Laplacian systems is constructed, and an iterative solve phase using multigrid cycles. General graphs pose algorithmic challenges not encountered in traditional multigrid applications. LAMG combines a lean piecewise-constant interpolation, judicious node aggregation based on a new node proximity measure (the affinity), and an energy correction of coarse-level systems. This results in fast convergence and substantial setup and memory savings. A serial LAMG implementation scaled linearly for a diverse set of \numgraphs real-world graphs with up to \maxedges edges, with no parameter tuning. LAMG was more robust than the UMFPACK direct solver and Combinatorial Multigrid (CMG), although CMG was faster than LAMG on average. Our methodology is extensible to eigenproblems and other graph computations.
\end{abstract}



\pagestyle{myheadings}
\thispagestyle{plain}
\markboth{O. E. LIVNE AND A. BRANDT}{LEAN ALGEBRAIC MULTIGRID FOR THE GRAPH LAPLACIAN}

\section{Introduction}
\label{introduction}
Let $G=(\cN,\cE,w)$ be a connected weighted undirected graph, where $\cN$ is a set of $n$ nodes, $\cE$ is a set of $m$ edges, and $w: \cE \rightarrow \Real^{+}$ is a weight function. The Laplacian matrix $\bA_{n \times n}$ is naturally defined by the {\it quadratic energy}
\be
    E(\bx) := \bx^{T} \bA \bx = \sum_{(u,v) \in \cE} w_{uv} \lp x_u-x_v \rp^2\,,\qquad \bx \in \Real^{\cN}\,,
    \label{energy}
\ee
where $\bx^T$ denotes the transpose of $\bx$. In matrix form,
\be
    \bA = \lp a_{uv} \rp_{u,v}\,,\quad
    a_{uv} :=
    \begin{cases}
        \sum_{v' \in \cE_u} w_{uv'}\,,& u=v\,, \\
        -w_{uv}\,,&v \in \cE_u := \left\{v': (u,v') \in \cE\right\} \,, \\
        0\,,& \text{otherwise.}
    \end{cases}
    \label{lap}
\ee
\begin{figure}[ht]
    \centering
    \begin{minipage}{.4\linewidth}
         \includegraphics[height=1in]{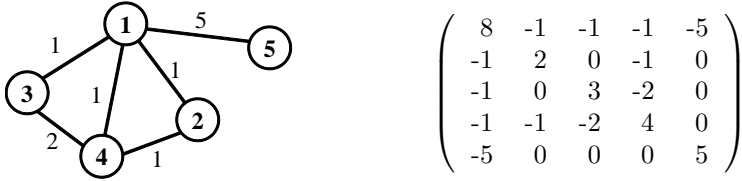}
    \end{minipage}
    \quad
    \begin{minipage}{.4\linewidth}
        $\left(
            \begin{tabular}{rrrrr}
            8 & -1 & -1 & -1  & -5 \\
           -1 &  2 &  0 & -1  &  0 \\
           -1 &  0 &  3 & -2  &  0 \\
           -1 & -1 & -2 &  4  &  0 \\
           -5 &  0 &  0 &  0  &  5
        \end{tabular}
        \right)$
    \end{minipage}
    \caption{A $5$-node graph and its corresponding Laplacian matrix.}
    \label{simple_graph}
\end{figure}

$\bA$ is Symmetric Positive Semi-definite (SPS), and has zero row sums and $2m+n$ non-zeros. Typically, $m \ll n^2$ and $\bA$ is sparse. Our approach also handles some SPS Laplacian matrices corresponding to {\it negative edge weights}, such as high-order and anisotropic grid discretizations \cite{alg_distance_anis, max_weight_basis}; those are discussed in \S \ref{negative_weights}.

Since $G$ is connected, $\bA$'s null space is spanned by the vector of ones $\bu$ (a disconnected graph can be decomposed into its components in $O(m)$ time \cite{sedgewick, tarjan}). We consider the nonsingular compatible linear system \cite[pp.~185--186]{trot}
\begin{subequations} \label{linsys}
\begin{eqnarray}
    &&\bA \bx = \bb \label{system} \\
    &&\bu^{T} \bx = 0\,, \label{xcompat}
\end{eqnarray}
\end{subequations}
where $\bb \in \Real^{\cN}$ is a given zero-sum vector, and $\bx \in \Real^{\cN}$ is the vector of unknowns.

Our goal is to develop an iterative numerical solver of (\ref{linsys}) that {\it requires $O(m)$ storage and $O(m \loget)$ operations to generate an $\ep$-accurate solution, for graphs arising in real-world applications.} The hidden constants should be small (in the hundreds, not millions). The solver should require a smaller cost to re-solve the system for multiple $\bb$'s -- a useful feature for time-dependent and other applications.

Importantly, we are interested in {\it good empirical performance} (bounded hidden constants over a diverse set of test instances), and do not consider the problem of designing an algorithm with provably linear complexity in any graph \cite[Problem 5, p. 18]{icm10}. While proofs are important, they often provide unrealistic bounds or no bounds at all on the hidden constants that can be attained in practice.

\subsection{Applications}
The linear system (\ref{linsys}) is fundamental to many applications; see Spielman's review \cite[\S 2]{icm10} for more details:
\bi
\item Elliptic Partial Differential Equations (PDEs) discretized on unstructured grids by finite elements within a fluid dynamics simulation \cite{boman, fischer}.
\item Interior-point methods for network flow linear programming \cite{DS08, FG07}.
\item Electrical flow through a resistor network $G$.
\ei
Additionally, $\bA \bx = \bb$ is {\it a stepping stone toward the eigenproblem.} Our multilevel methodology can be extended to compute the smallest eigenpairs of $\bA$ with minor adaptations (cf.~\S \ref{eigenvalue}). The Laplacian eigenproblem is central to graph regression and classification in machine learning \cite{hongkong_class, zhu_learning}, spectral clustering of images, graph embedding \cite{radu}, and dimension reduction for genetic ancestry discovery \cite{ancestry}. Of particular interest is the Fiedler value -- the smallest non-zero eigenvalue of $\bA$, which measures the algebraic connectivity of $G$ \cite[\S 1.1]{chung} and is related to minimum cuts \cite{ding}. Although we believe it is preferable to develop multiscale strategies for the original formulations of these problems, as demonstrated by the works \cite{class04, clustering06, alg_distance, nature} and graph partitioning packages \cite{chaco}, a fast black-box eigensolver is a practical alternative.

\subsection{Related Work}
There are two main approaches to solving (\ref{linsys}): direct, leading to an exact solution (up to round-off errors); and iterative, which typically requires a one-time setup cost, followed by a solve phase that produces successive approximations $\tbx$ to $\bx$ to achieve $\ep$-accuracy, namely,
\be
	\|\bx-\tbx\|_{\bA} \leq \ep \|\bx\|_{\bA}\,,\qquad \|\bx\|_{\bA} := \sqrt{E(\bx)}\,.
	\label{epaccuracy}
\ee

\subsubsection{Direct Methods}
The Cholesky factorization with a clever elimination order can be applied to $\bA$. A permutation matrix $\bP$ is chosen and factorization $\bP^{T} \bA \bP = \bL \bL^{T}$ constructed so that the lower triangular $\bL$ is as sparse as possible, using Minimum or Approximate Minimum Degree Ordering \cite{add, md}. Except for simple graphs, direct algorithms do not scale, requiring $O(n^{1.5})$ operations for planar graphs \cite{nested_dissection, generalized_nested_dissection} and $O(n^3)$ in general. Alternatively, fast matrix inversion can be performed in $O(n^{2.376})$ or combined with Cholesky, yet yields similar complexities \cite[\S 3.1]{icm10}.

\subsubsection{Iterative Methods: Graph Theoretic}
These are variants of the preconditioned conjugate method \cite[\S 10.3]{gvl} that achieve (\ref{epaccuracy}) in $O(\sqrt{\kappa(\bA \bB^{-1})}$ $\loget)$ iterations for a preconditioner $\bB$; $\kappa$ is the finite condition number \cite[\S 3.3]{icm10}.

Spielman and Teng (S-T) \cite{st06} and subsequent works \cite{koutis} have been focusing on multilevel graph-sparsifying preconditioners. The S-T setup builds increasingly smaller graphs, alternating between partial Cholesky and ultra-sparsification steps, in which the graph is partitioned into sets of high-conductance nodes without removing too many edges. The complexity is a near-linear $O(m \log^2 n \log(1/\ep))$, guaranteed for any symmetric diagonally-dominant $\bA$. Unfortunately, no implementation is available yet, nor is there a guarantee on the size of the hidden constant.

\subsubsection{Iterative Methods: AMG}
Algebraic Multigrid (AMG) is a class of high-performance linear solvers, originating in the early 1980s \cite[\S 1.1]{guide}, \cite{geodetic}, \cite{rs86} and under active development. During setup, AMG recursively constructs a multi-level hierarchy of increasingly coarser graphs by examining matrix entries, without relying on geometric information. The solve phase consists of multigrid cycles. AMG can be employed either as a solver or a preconditioner \cite[App.~A]{trot}. Open-source parallel implementations include Hypre \cite{hypre} and Trilinos-PETSc \cite{trilinos}. In classical AMG, the coarse set is a subset of $\cN$; alternatively aggregation AMG \cite[App.~A.9]{trot},\cite{braess, blatt, notay_etna} and smoothed aggregation \cite{sa} define coarse nodes as aggregates of fine nodes.

AMG mainly targets discretized PDEs, where it has been successful \cite{fischer}. Advanced techniques have ventured to widen its scope by increasing the interpolation accuracy. These include Bootstrap AMG \cite[\S 17.2]{yes}, \cite{bamg} adaptive smoothed aggregation \cite{asa}, and interpolation energy minimization \cite{olson}. While these methods approach linear scaling for more systems, their complexity cannot be controlled in general graphs (cf. \S \ref{interpolation_caliber}).

At the same time, accelerated aggregation AMG has become a hot research topic. A crude {\it caliber-1} (piecewise-constant) interpolation is employed between levels to reduce runtime and memory costs, at the expense of a slower cycle that is subsequently accelerated. Notay \cite{notay_etna, notay_nonsym} aggregated nodes based on matrix entries and applied multilevel CG acceleration with a large cycle index to obtain a near-optimal solver for convection-diffusion M-matrices, but the method was limited to those grid graphs. Caliber-1 interpolation was tested for a single graph by Bolten et. al within the bootstrap framework \cite{bamg_markov}, but their setup cost was large and required parameter tuning.

The present work aims at generalizing AMG to graph Laplacians and addresses peculiarities not encountered in traditional AMG applications. To the best of our knowledge, the only other solver targeting general topologies is Combinatorial Multigrid (CMG) \cite{cmg}, a hybrid graph-theoretic-AMG preconditioner that partitions nodes into high-conductance aggregates, similarly to S-T. CMG outperformed classical AMG for a set of 3-D image segmentation applications. In our experiments over a much larger graph collection, CMG and LAMG had comparable average solve speeds, yet LAMG's performance was much more robust, with almost no outliers.

\subsection{Our Contribution}
We present Lean Algebraic Multigrid (LAMG): a practical graph Laplacian linear solver. A \matlab LAMG implementation scaled linearly for a set of \numgraphs real-world graphs with up to \maxedges edges, ranging from computational fluid dynamics to web, biological and social networks. Specifically, the setup phase required on average $\setupmvm$ Matrix-Vector Multiplications (MVMs) and $\storageperedge m$ storage bytes, and the average solve time was $\approx \solvemvm \log(1/\ep)$ MVMs per right-hand side. The standard deviations were small with only three outliers. LAMG was more robust than the UMFPACK direct solver and CMG, although CMG was faster on average (\S \ref{smorgasbord}). Our methodology is extensible beyond the scope of S-T and CMG, to non-diagonally-dominant (\S \ref{negative_weights}), eigenvalue, and nonlinear problems (\S \ref{eigenvalue}).

LAMG is an accelerated caliber-1 aggregation-AMG algorithm that builds upon the state-of-the-art AMG variants, yet introduces four new ideas essential to attaining optimal efficiency in general graphs:
\bi
    \item[(a)] {\it Lean methodology} (\S \ref{lean}). LAMG advocates using {\it minimalistic over sophisticated AMG components}. No parameter tuning should be required. In particular, we apply caliber-1 interpolation between levels, constructed using relaxed Test Vectors (TVs), but without bootstrapping them as in the papers \cite{bamg_markov, bamg}. Fast asymptotic convergence is achieved by the following three ideas.
    \item[(b)] {\it Low-degree Elimination} (\S \ref{elimination}). Like the S-T method \cite{st06}, we eliminate low-degree nodes prior to each aggregation. However, the role of elimination here is different: {\it it removes the effectively-1-D part of the graph}, thereby eliminating extreme tradeoffs between complexity and accuracy in aggregating the remaining graph. Thus, unlike S-T, our elimination need not strictly reduce the number of edges (allowing us to eliminate nodes of higher degrees) nor be exact (making it also useful in the eigenproblem). The elimination and aggregation are in fact specializations of the same coarsening scheme (\S \ref{eigenvalue}).
    \item[(c)] {\it Affinity} (\S \ref{affinity}). Aggregation is based on a new {\it normalized relaxation-based node proximity} heuristic. Ron et al. \cite{alg_distance} were the first to use TVs to measure algebraic distance between nodes, but our measure is effective for a wider variety of graph structures. The affinity admits a statistical interpretation and approximates the diffusion distance \cite{diffusion_maps}. In contrast, the S-T algorithm strives to create high-conductance aggregates \cite{st06, koutis_conductance}.
    \item[(d)] {\it Energy-corrected Aggregation} (\S \ref{aggregation}). Recognizing that caliber-1 interpolation leads to an energy inflation of the coarse-level system (\S \ref{inflation}), our aggregation also {\it reuses test vectors to minimize coarse-to-fine energy ratios}. We offer two alternative energy corrections to accelerate the solution cycle: a flat correction to the Galerkin operator resemblant of Braess' work \cite{braess}, yet resulting in a superior efficiency; and an adaptive correction to the solution via {\it multilevel iterate recombination} \cite[\S 7.8.2]{trot} that is even more efficient.
\ei

Following an AMG prelude in \S \ref{basics}, we discuss each of the main ideas in \S \ref{ideas}. They are integrated into the complete LAMG algorithm in \S \ref{algorithm} (cf. the expanded ArXiV e-print \cite{lamg_arxiv} for implementation details). Our development methodology emphasizes learning from examples: we studied instances for which the original design was slow to derive general aggregation rules. Testing over a large collection ensured that new rules did not spoil previous successes. The results are presented in \S \ref{results}. Coarsening improvements and extensions to the eigenproblem and other graph computational problems are outlined in \S \ref{extensions}.

\section{Algebraic Multigrid Basics}
\label{basics}
Relaxation methods for $\bA\bx=\bb$ such as Jacobi or Gauss-Seidel slow down asymptotically. Yet after only several sweeps, the {\it error} $\bee := \bx-\tbx$ in the approximation $\tbx$ to $\bx$ becomes {\it algebraically smooth}: its normalized residuals are much smaller than its magnitude \cite[\S 1.1]{guide} (assuming its mean has been subtracted out). In AMG, these errors are approximated by an interpolation $\bP_{n \times n_c}$ from a coarse subspace: $\bee \approx \bP \bee^c$, where $\bee^c$ is a coarse vector of size $n_c$.

Recognizing that $\bA\bx=\bb$ corresponds to the quadratic minimization
\be
    \bx = \mymin{\by} \Etot(\by)\,,\quad \Etot(\bx) := \frac12 \bx^{T} \bA \bx - \bx^{T} \bb\,,
    \label{quad}
\ee
the variational correction scheme \cite{geodetic} seeks the optimal correction in the energy norm,
\be
    \bee^c = \argmin{\by^c} \Etot\left(\tbx + \bP \by^c \right)\,.
    \label{min_ec}
\ee
The resulting two-level cycle (except for determining $\bx$'s mean) is
\ben
    \item Perform 2-3 relaxation sweeps on $\bA\bx=\bb$, resulting in $\tbx$.
    \item Compute an approximation $\tbee^c$ to the solution $\bee^c$ of
    \be
        \bA^c \bee^c = \bb^c\,,\quad \bA^c := \bP^{T} \bA \bP \,,\quad \bb^c := \bP^{T} \left(\bb - \bA \tbx \right)\,.
        \label{galerkin}
    \ee
    \item Correct the fine-level solution:
    \be
        \tbx \leftarrow \tbx + \bP \tbee^c\,.
        \label{correction}
    \ee
\een
Eq.~(\ref{galerkin}), called Galerkin coarsening, is a smaller linear system. In the multilevel cycle, it is recursively solved using $\gamma$ two-level cycles, where the {\it cycle index} $\gamma$ is an input parameter. Our interpolation $\bP$ is full-rank with unit row sums (see \S \ref{exten}) for which it easy to verify that $\bA^c$ is also a connected graph Laplacian. Overall, LAMG constructs a hierarchy of $L$ increasingly coarser Laplacian systems (``levels'') $\bA^l \bx^l=\bb^l$, $l=1,\dots,L$, the finest being the original system $\bA^1:=\bA, \bb^1:=\bb$. The setup phase depends on $\bA$ only, and produces $\{(\bA^l,\bP^l)\}_{l=2}^L$, where $\bA^l$ is $n_l \times n_l$ and $\bP^l$ is the $n_{l-1} \times n_l$ interpolation matrix from level $l$ to $l-1$. We denote by $G^l=(\cN^l,\cE^l)$ the graph corresponding to $\bA^l$.

\section{LAMG: Main Ideas}
\label{ideas}
Our description refers to a single coarsening stage ($\bA \rightarrow \bA^c$) and applies to each pair of levels.

\subsection{Lean Methodology}
\label{lean}

\subsubsection{Relaxation}
\label{relax}
Our choice is Gauss-Seidel (GS) relaxation, defined by the successive updates \cite[\S 1.1]{guide}
\be
    \text{For } u=1,\dots,n\,,\qquad x_u \leftarrow \frac{b_u - \sum_{v \in \cE_u} a_{uv} x_v}{a_{uu}}\,.
    \label{gs}
\ee
We picked GS because it is an effective smoother in SPS systems \cite[\S 1]{guide} and does not require parameter tuning (such as the Jacobi relaxation damping parameter \cite{bamg_markov}).

\subsubsection{Interpolation Caliber}
\label{interpolation_caliber}
Textbook multigrid convergence for the Poisson equation requires that the interpolation of corrections $\bP$ be second-order \cite[\S 3.3]{mg_theory}. The analogous AMG theory implies a similar condition on the interpolation accuracy of low-energy errors. While a piecewise-constant $\bP$ is acceptable in a two-level cycle, it is insufficient for V-cycles ($\gamma=1$); W-cycles ($\gamma=2$) are faster but costly \cite[p.~471]{trot},\cite{notay_etna}.

Constructing a second-order $\bP$ is already challenging in grid graphs \cite{bamg, sa, olson}. We argue that it is infeasible in general graphs:
\ben
    \item[(a)] The graph's effective dimension $d$ is unknown; had it been known, the required interpolation {\it caliber}, i.e., the number of coarse nodes used to interpolate a fine node, would grow with $d$ and result in unbounded interpolation complexity.
    \item[(b)] Identifying a proper interpolation set (whose ``convex hull'' contains the fine node) is a complex and costly process \cite{bamg}.
    \item[(c)] The Galerkin coarse operator $\bP^{T} \bA \bP$ fills in considerably; cf. \S\ref{fillin}.
\een
In contrast, LAMG employs a lean caliber-1 (piecewise-constant) $\bP$, equivalent to an {\it aggregation} of the nodes into coarse-level aggregates, and corrects the energy of the coarse-level Galerkin {\it operator} to maintain good convergence (in practice, the correction is applied to the coarse right-hand side). This could not have been achieved within the variational setting of \S \ref{basics}, which only permits modifying $\bP$. Here the barrier to fast convergence is the coarse-to-fine operator energy ratio. Our contribution is an algorithm that yields a small energy ratio, which translates into optimal efficiency; cf.~\S \ref{aggregation}. The compatible relaxation performance predictor \cite{cr_etna, cr_oren, cr_james}, \cite[\S\S 14.2--14.3]{guide} is irrelevant for low interpolation accuracy; the energy ratio is a better predictor.

\subsubsection{Managing Fill-in}
\label{fillin}
Frequently, the coarse-level matrices in AMG hierarchies become increasingly dense. This is a result of a poor aggregation, a high-caliber $\bP$, or both: many fine nodes whose neighbor sets are disjoint are aggregated, creating additional edges among coarse-level aggregates. This renders the ideally-accurate interpolation irrelevant, because the actual cycle {\it efficiency} (error reduction per unit work) is small albeit convergence may be rapid. While fill-in is often controllable in grid graphs because their coarsening is still local, it is detrimental in non-local graphs.

LAMG's interpolation is designed to minimize fill-in. Heuristically, the sparser $\bP$, the sparser $\bP^{T} \bA \bP$. We further indirectly control fill-in via our affinity criterion (\S \ref{affinity}), which tends to aggregate nodes that share many neighbors. The cycle work is also restrained by a fractional cycle index \cite[\S 6.2]{guide} between 1 and 2; cf.~\S \S \ref{setup},\ref{solve}.

Occasionally, the interpolation caliber may be slightly increased as long as the number of coarse edges does not become too large; see \S \ref{improvements}.

\subsubsection{Utilizing Extenuating Circumstances}
\label{exten}
Specific properties of graph Laplacians are exploited to simplify the LAMG construction.
\bi
    \item Since $\bA$ has zero row sums, its null-space consists of constant vectors. A fundamental AMG assumption is that {\it all near}-null-space errors can be fitted by a {\it single} interpolation from a coarse level \cite[p.~8]{guide}. Here it implies that the interpolation weights can be apriori set to $1$. The unit-weight assumption is easily verified for Laplacians with bounded node degrees \cite[p.~439]{trot}; in the most interesting applications of this work, however, the node degree is unbounded, and a proof of this conjecture is an open problem.
    \item Some graph locales are effectively one-dimensional: many nodes have degree 1--2. Such nodes can be quickly eliminated similarly to the paper \cite{koutis} (\S \ref{elimination}).
    \item In other graphs, Gauss-Seidel is an efficient solver and no coarsening is required. These include complete graphs, star graphs, expander graphs \cite[\S 1]{icm10}, and certain classes of random graphs. More generally, if GS converges fast for a subset of the nodes, they can all be aggregated together (\S \ref{agg_algorithm}).
\ei

\subsection{Low-degree Elimination}
\label{elimination}
We first attempt to eliminate from $\cN$ an independent set $\cF$ of nodes $u$ of degree $|\cE_u| \leq 4$. The set is identified by initially marking all nodes as ``eligible''; we then sweep through nodes, adding each eligible low-degree node to $\cF$ and marking its neighbors as ineligible \cite[Algorithm~1]{lamg_arxiv}.

Eliminating a node connects all its neighbors; therefore, $\cF$-nodes of degree $\leq 3$ do not increase $m$. When $|\cE_u| = 4$, $m$ might be increased by at most $2$. However, we assume that this is unlikely to happen for many $\cF$-nodes and eliminate these nodes as well (in practice, we have observed that the neighbors are already connected prior to elimination). Eliminating larger degrees results in an impractical fill-in.

Let $\cC := \cN \backslash \cF$. $\bA \bx = \bb$ reduces to the Schur complement system
\be
    \bA^c \bx_{\cC} = \bb^c\,,\,\,\,\, \bA^c := \bP^{T} \bA \bP\,,\,\,\,\, \bb^c := \bP^{T} \bb\,,\,\,\,\,
    \bP := \bPi \left(-\bA_{\cF\cC}^{T} \bA_{\cF\cF}^{-1}\,,\bI_{\cC} \right)^{T}\,,
    \label{elim}
\ee
where $\bPi$ is a permutation matrix such that $\bPi^{T} \bx$ lists the all $\cF$ node values, then all $\cC$ node values. (\ref{elim}) is a smaller Laplacian system for which we perform further elimination rounds, until $|\cF|$ becomes small \cite[Algorithm~2]{lamg_arxiv}.

The purpose of elimination is to reduce $n$ while incurring a small fill-in, and to remove the 1-D part of the graph, which cannot be effectively coarsened by the energy-corrected aggregation of \S \ref{aggregation}. Note that $\bP$'s caliber is larger than $1$ here.

Viewed as a full approximation scheme (\S \ref{eigenvalue}), (\ref{elim}) could be generalized to a {\it non-exact} elimination for approximating the lowest eigenvectors of $\bA$, instead of forming a nonlinear Schur complement \cite{amls}. Additionally, a larger set of {\it loosely-coupled} nodes could also be eliminated: if $\bA_{\cF\cF}$ is strongly diagonally-dominant, its inverse can be approximated by a few Jacobi relaxations. The two-level convergence rate and fill-in would need to be kept in check in that case. We plan to pursue these generalizations in a future research.

We hereafter denote the coarse system by $\bA \bx = \bb$ (either (\ref{system}), if $q=0$, or (\ref{elim})), which is further coarsened by caliber-1 aggregation in \S\S \ref{affinity}--\ref{aggregation}.

\subsection{Affinity}
\label{affinity}

The construction of an effective aggregate set hinges upon defining which nodes in $\cN$ are ``proximal'', i.e., nodes whose values are strongly coupled in all smooth (i.e., low-energy) vectors \cite[p.~473]{trot}. Table~\ref{prox_def} lists three definitions.
\begin{table}[htbp]
\centering
\begin{tabular}{|l|l|}
    \hline
    Classical AMG \cite{rs86} & $1 - |w_{uv}| / \max\left\{ \max_s |w_{us}|, \max_s |w_{sv}| \right\}$ \\ \hline
    Algebraic Distance \cite{alg_distance} & $ \max_{k=1}^K \left| x_u^{(k)} - x_v^{(k)} \right|$ \\ \hline
    \multirow{2}{*}{Affinity (LAMG)} & $c_{uv} := 1 - \left|\left(X_u,X_v \right)\right|^2/\left(\left(X_u,X_u \right) \left(X_v,X_v \right)\right)\,,$ \\
    & $(X_u,X_v):= \sum_{k=1}^K x^{(k)}_u x^{(k)}_v$ \\
    \hline
\end{tabular}
\caption{Comparison of node proximity measures. Nodes are defined as ``close'' when the measure is smaller than a threshold. $\{\bx^{(k)}\}_{k=1}^K$ is a set of relaxed test vectors; see the text.}
\label{prox_def}
\end{table}
\begin{figure}[htbp]
    \centering
    \begin{center}
    \begin{tabular}{cc}
      \includegraphics[height=1.25in]{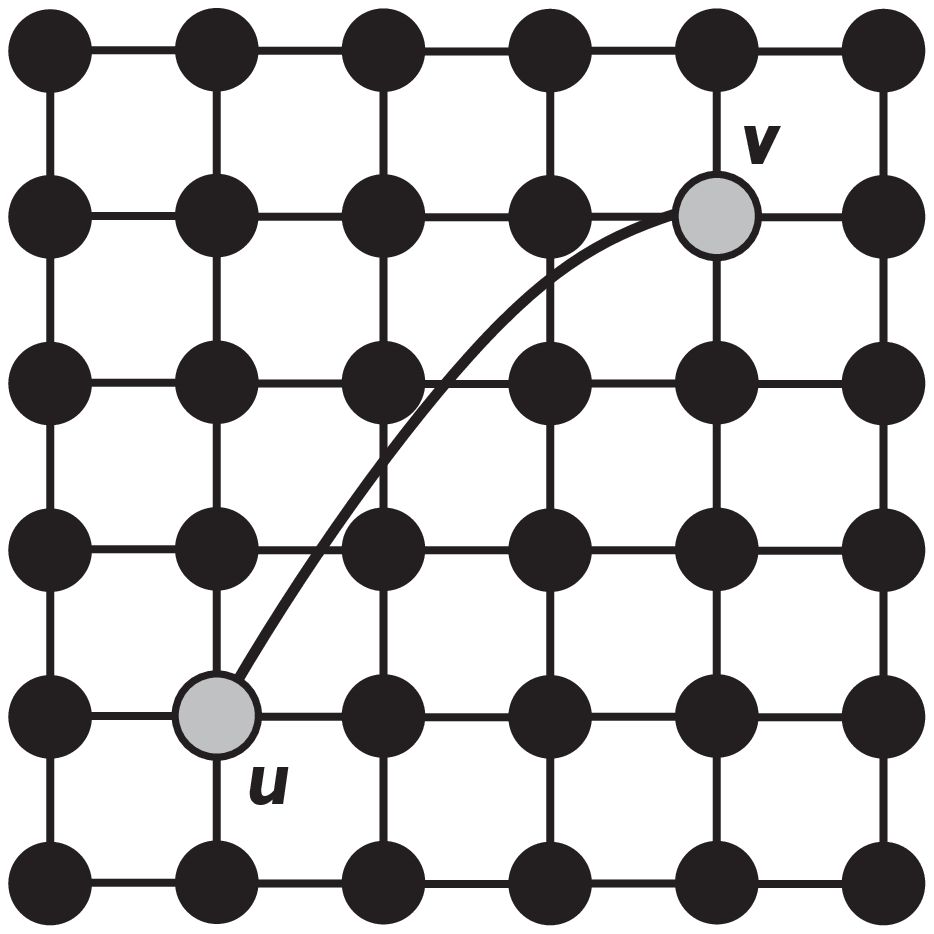}
      &
      \includegraphics[height=1.2in]{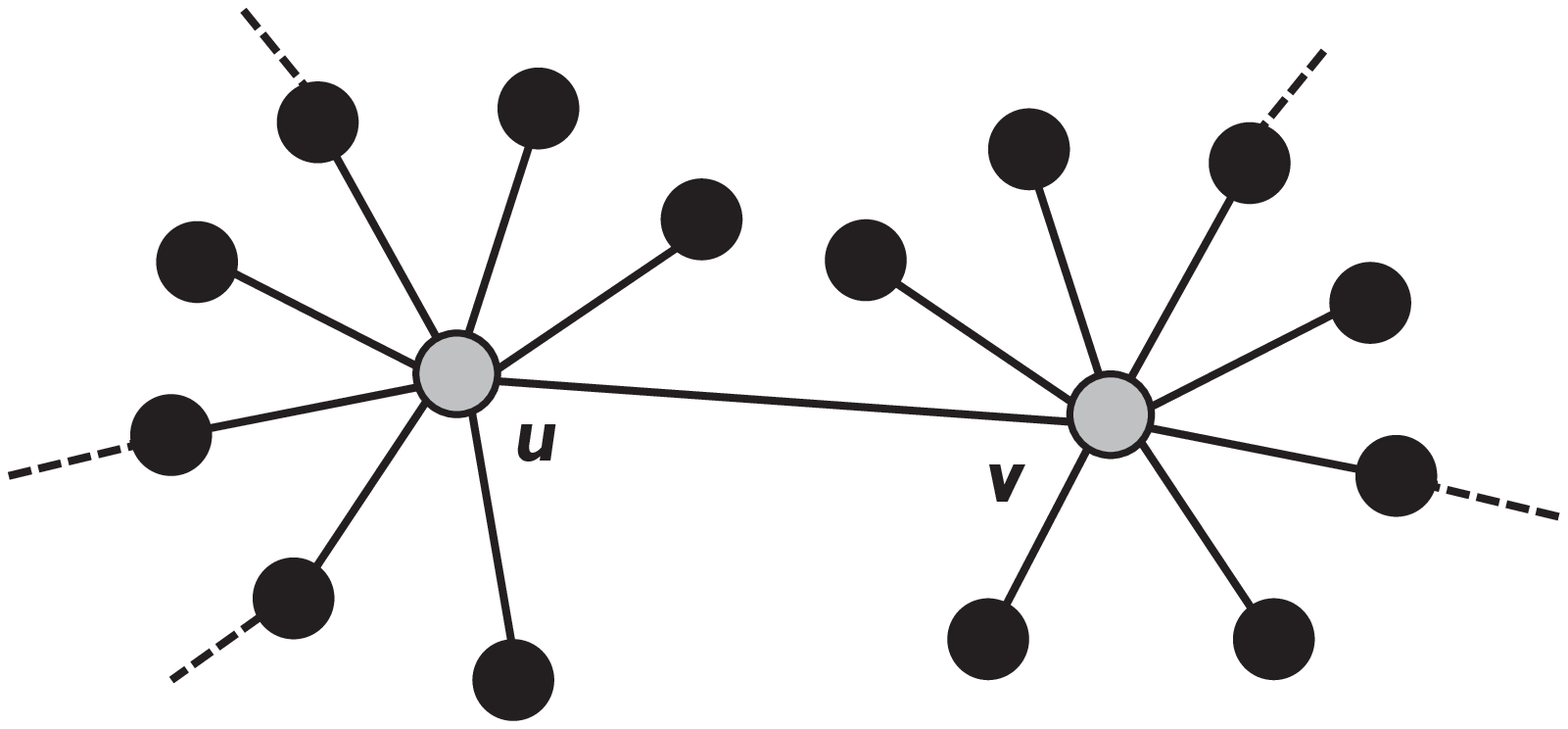}
      \\
      (a) & (b) \\
    \end{tabular}
    \end{center}
    \caption{Unweighted graph instances that present aggregation difficulties. (a) A 2-D grid with an extra link. (b) Two connected hubs. A hub is a high-degree node.}
    \label{two_examples}
\end{figure}

\subsubsection{Existing Proximity Measures}
Classical AMG defines proximity based on edge weights (Table~\ref{two_examples}, top row). While this has worked well for coarsening discretized scalar elliptic PDEs, it leads to wrong aggregation decisions in non-local graphs. In a grid graph with an extra link between distant nodes $u$, $v$ (Fig.~\ref{two_examples}a), $u$ and $v$ become proximal and may be aggregated. Unless $w_{uv}$ is outstandingly large, this is undesirable because $u$ and $v$ belong to unrelated milieus of the grid.

This problem is overcome by the algebraic distance measure introduced by Ron et al. \cite{alg_distance} (Table~\ref{prox_def}, middle row; a related definition is used in the work \cite{alg_distance_anis}). Insofar as coarsening concerns the space of smooth error vectors $\bx$, nodes $u$ and $v$ should be aggregated only if $\bx_u$ and $\bx_v$ are highly correlated for all such $\bx$. A set of $K$ Test Vectors (TVs) $\bx^{(1)},\dots,\bx^{(K)}$ is generated -- a sample of this error space \cite[\S 17.2]{yes}, \cite{alg_distance}. Each TV is the result of applying $\nu$ relaxation sweeps to $\bA \bx = \bzero$, starting from $\mathrm{random}[-1,1]$. However, Ron et al.'s definition falls prey to a graph containing two connected high-degree nodes $u$ and $v$ (``hubs''; Fig.~\ref{two_examples}b). For each $k$, the value $x^{(k)}_u$ is an average over a large neighborhood of random node values whose size increases with the number of sweeps, hence is small. Similarly, $x^{(k)}_v$ is small, so every such $u$ and $v$ turns out proximal, even though they may be distant.

\subsubsection{The New Proximity Measure}
LAMG's proximity measure, the {\it affinity} (Table~\ref{prox_def}, bottom), also relies on TVs, but is scale-invariant, and correctly assesses both Fig.~\ref{two_examples}a and b as well as many other constellations. The affinity $c_{uv}$ between $u$ and $v$ is defined as the goodness of fitting the linear model $x_v \approx p\, x_u$ to TV values:
\be
    c_{uv} := 1 - \frac{\left|\left(X_u, X_v\right)\right|^2}{\left(X_u, X_u\right) \left(X_v, X_v\right)}\,,\,\,\,\,\,
    \left( X, Y \right) := \sum_{k=1}^K x^{(k)} y^{(k)}\,,\,\,\,\, X_u := \left(x^{(1)}_u,\dots,x^{(K)}_u\right)\,.
    \label{cuv}
\ee
$c_{uu} = 0$, $0 \leq c_{uv} \leq 1$ and $c_{uv}=c_{vu}$. The affinity measures {\it distance}: the smaller $c_{uv}$, the closer $u$ and $v$. In the $d$-D discretized Laplace operator on a grid, $c_{uv}$ related to the {\it geometric distance} between the gridpoints corresponding to $u$ and $v$. In general graphs, $c_{uv}$ is an alternative definition of the algebraic distance \cite{alg_distance}, and approximates the {\it diffusion distance} at a short time $\nu$ \cite{diffusion_maps}.

\subsubsection{Statistical Interpretation}
$x_u$ can be thought of as a random variable; $c_{uv}$ is the Fraction of Variance Unexplained of linearly regressing $x_u$ on $x_v$ using the TV samples $X_u$ and $X_v$ \cite{regression}. We make several observations:
\bi
    \item $c_{uv}$ is invariant to scaling $X_u$ and $X_v$. This is vital in the two-hub case (Fig.~\ref{two_examples}b).
    \item Rather than subtracting the sample means $\overline{X}_u := \sum_{k=1}^K x^{(k)}_u$ and $\overline{X}_v$ from $X_u$ and $X_v$, respectively, as in the standard statistical definition, we use their exact means over all error vectors. These are zero, since each $X_u$ is a linear combination of some initial $X_w$'s, each of which has a zero mean.
    \item A weighted inner product $(X_u,X_v)$ could be used to account for TV variances, but since all TVs have the same level of smoothness, i.e., comparable normalized residuals \cite[\S 17.2]{yes}, \cite{bamg}, they were assigned equal weights.
\ei
A few vectors $K$ and smoothing sweeps $\nu$ suffice to obtain a good enough $c_{uv}$ estimate, which guides a coarsening of the node set by a modest factor of 2--3 (cf.~\S \ref{agg_algorithm}).

\subsubsection{Interpolation Accuracy}
Bootstrap AMG \cite[\S 17.2]{yes} defines a general-caliber $\bP$ using a least-squares fit to TVs. In our case,
\be
    c_{uv} = \mymin{p} \frac{\left\|p X_u-X_v\right\|^2}{\left\|X_v\right\|^2}
    \label{interp_accuracy}
\ee
relates the affinity to the accuracy of the caliber-$1$ interpolation formula $x_v = p x_u$ for TVs. As a byproduct, we obtain the interpolation coefficient $\hat{p} = (X_u,X_v)/(X_v,X_v)$. Eq.~(\ref{cuv}) also works for non-zero-sum matrices, such as restricted and normalized Laplacians \cite[\S 2]{icm10}, where $p_{uv}$ is set to $\hat{p}$ (see also \S\S \ref{improvements},\ref{other}). For the Laplacian, $\hat{p}$ is abandoned in favor of $p_{uv}=1$ (cf.~\S \ref{exten}).

In the Helmholtz equation, $c_{uv}$ is large for all $u,v$, indicating that all nodes are distant and that no single aggregate set can yield fast AMG convergence (indeed, multiple coarse grids are required to restore textbook multigrid efficiency \cite{wave_accuracy}).

\subsection{Aggregation}
\label{aggregation}
Aggregation levels refer to the two-level method of \S \ref{basics} with a caliber-$1$ interpolation. $\bP$  is equivalent to partitioning $\cN$ into $n_c$ non-overlapping {\it aggregates} $\{\cT_U\}_{U \in \cN^c}$; $\cT_U$ is the set of $\bee_u$ interpolated from $\bee^c_U$, and $\cN^c := \{1,\dots,n_c\}$. Each aggregate consists of a {\it seed} node and zero or more {\it associate} nodes \cite[Fig.~3.2]{lamg_arxiv}.

This section explains the technical details of aggregate selection, and may be of interest to multigrid experts. Other readers may wish to skip it and assess the quality of our aggregation decisions via the numerical experiments in \S \ref{results}.

\subsubsection{Aggregation Rules}
Intuitively, nodes should be aggregated together if their values are ``close''; ideal aggregates have strong internal affinities and weaker external affinities. To this end, we formulated five rules:
\ben
    \item Each node can be associated with one seed.
    \item A seed cannot be associated.
    \item Aggregate together nodes with smaller affinities before larger affinities.
    \item Favor aggregates with small energy ratios (\S \ref{energy_correction}).
    \item A hub node should be a seed.
\een
Rules 1 and 2 prevent an associate from being transitively with a distant seed. Otherwise, long chains might be aggregated together, creating aggregates with weak internal connections and very large energy ratios. Rule 3 favors strongly-connected aggregates. Rule 4 has a dual purpose: (a) Ultimately, the energy ratio determines the AMG asymptotic convergence factor; hence, this rule ensures good convergence. (b) Affinities are based on local information (relaxed TVs); their quantitative value becomes fuzzier as nodes grow apart \cite[\S 5]{alg_distance}. Since small energy ratios usually dictate small aggregates, affinities are thus used only for {\it local} aggregation decisions. Rule 5 avoids costly aggregation decisions due to traversing the large neighbor sets of a hubs, which would increase the computational work (cf.~\S \ref{energy_correction}). On the other hand, aggressively coarsening a large clique into a single node is desirable, as all its nodes are strongly correlated. Thus hubs are defined as {\it locally-high degree} nodes.\footnote{In our code, a hub is a node whose degree is significantly larger than a weighted mean of its neighbors' degrees \cite{communities}: $|\cE_u| \geq 8 \sum_{v \in \cE_u} |w_{uv}| |\cE_v|/\sum_{v \in \cE_u} |w_{uv}|$.}

A typical coarsening ratio in our algorithm ranges between .3--.5.

\subsubsection{Aggregation Algorithm}
\label{agg_algorithm}
The algorithm requires the cycle index $\gamma$ as an input. Each node is marked as a seed, associate or undecided. First, hubs are identified and marked as seeds. Second, edges with very small $|w_{uv}|$ are discarded during aggregation. Since relaxation converges fast at the nodes that become disconnected, they need not be coarsened at all; however, to keep the coarse-level matrix a Laplacian, we aggregate all of them into a single (dummy) aggregate. All other nodes are marked as undecided.

Aggregation is performed in $r$ stages: aggregate sets $S_1,\dots,S_r$ are generated such that each $S_i$-aggregate is contained in some $S_{i+1}$-aggregate. The set whose coarsening ratio $\alpha:=|S_i|/n$ is closest to $\amax := .7/\gamma$ is selected as the final set, so that the total cycle work would be bounded by $\approx 1 + \gamma \amax + (\gamma \amax)^2 + \dots \approx \frac{10}{3}$ finest-level units, had the same fill-in occurred at all levels (this is a practical guideline for selecting a good coarse set; there exist many sets that yield similar cycle complexities). In our code, at most two stages are performed per coarsening level in the cycle.

In each stage, we scan undecided nodes $u$ and decide whether to aggregate each one with a neighbor $s$ that is either an existing seed, or an undecided node that thereby becomes a new seed. $s$ is the non-associate neighbor of $u$ with the smallest affinity $c_{us}$. At the end of the last stage, still-undecided nodes are converted to seeds. The complete algorithm is described in the ArXiV e-print \cite[\S 3.4.3]{lamg_arxiv}.

\subsubsection{Energy Inflation}
\label{inflation}
\begin{figure}[htbp]
    \centering
    \begin{center}
    \begin{tabular}{cc}
      \includegraphics[height=1.5in]{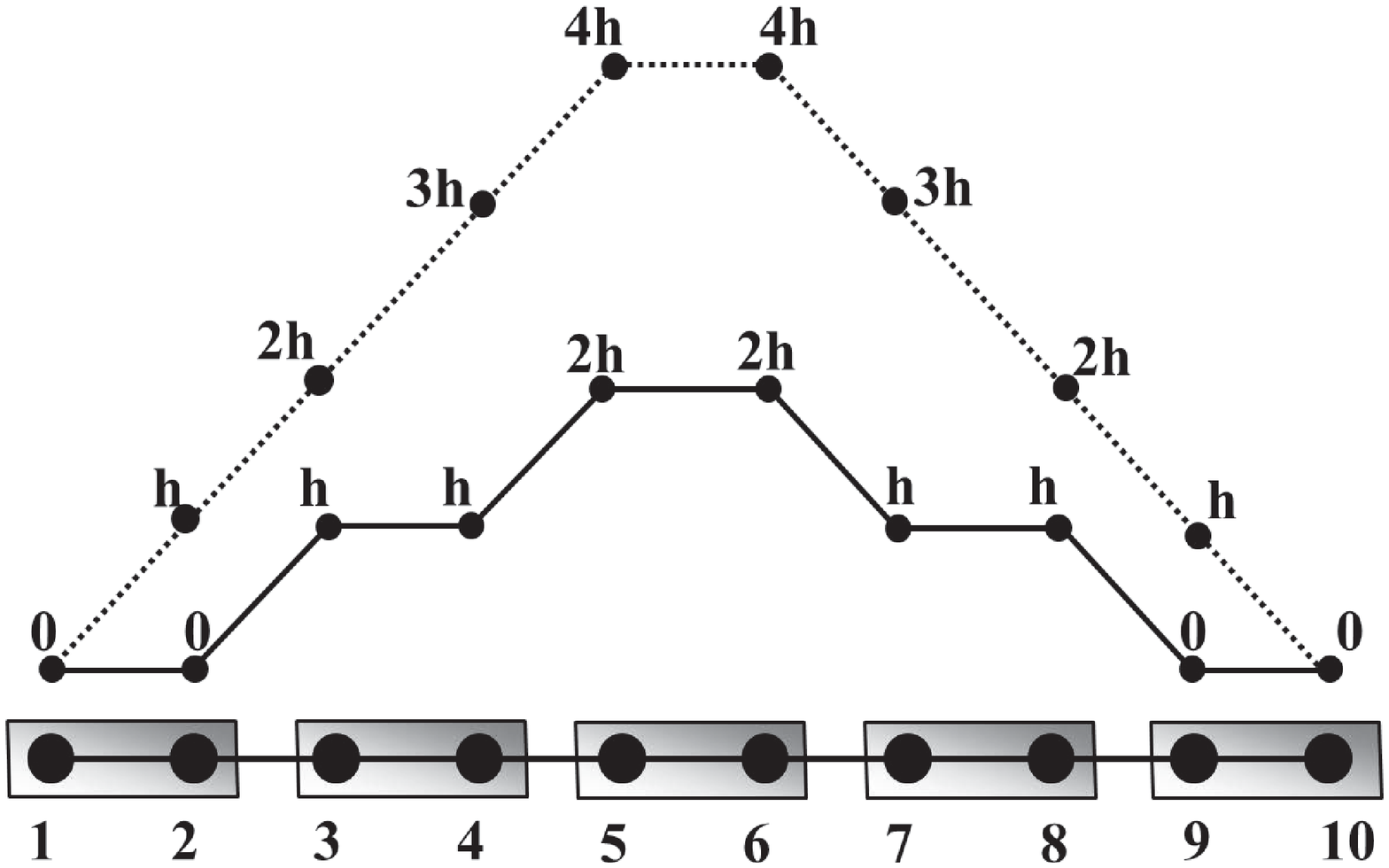}
      &
      \includegraphics[height=1.5in]{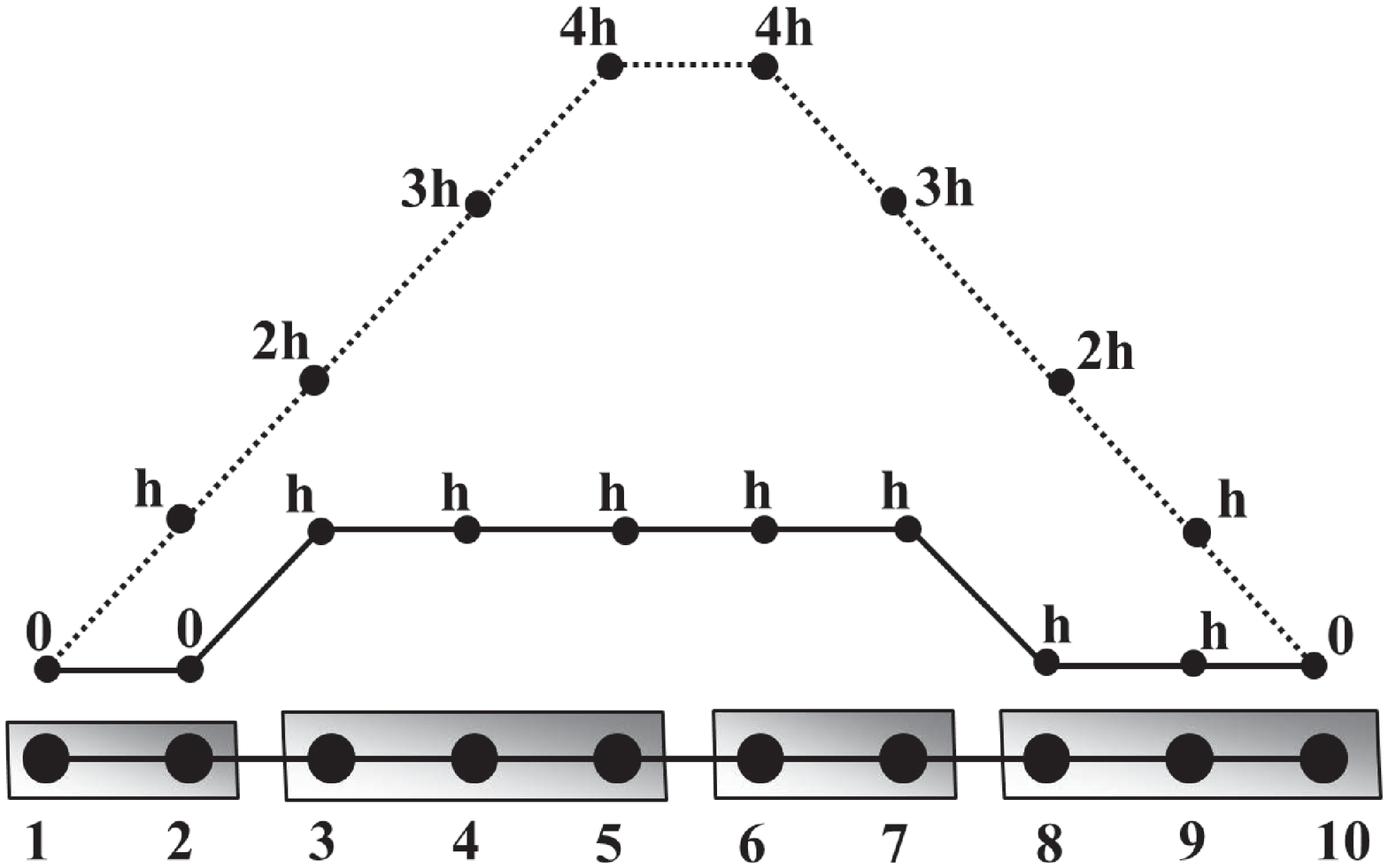}
      \\
      (a) & (b) \\
    \end{tabular}
    \end{center}
    \caption{The Galerkin correction $P \bee^c$ to a piecewise linear error $\bee$ in a path graph for (a) uniform 1:2 aggregation and (b) variable aggregate size. The meshsize $h>0$ is arbitrary.}
    \label{path}
\end{figure}

The Galerkin coarse-level correction $\bee^c$ (Eq.~(\ref{min_ec})) is the best approximation to a smooth error $\bee$ in the energy norm. Braess \cite{braess} noted that this does not guarantee a good approximation in the $l_2$ norm. For example, if $\bee$ is a piecewise linear function in a path graph (1-D grid with $w \equiv 1$) coarsened by aggregates of size two, $\bP \bee^c$ is constant on each aggregate and matches $\bee$'s slope across aggregates, resulting in about half the fine-level magnitude. See Fig.~\ref{path}a.

An equivalent and more useful observation is that the energy of $\bP \bT \bee$ is twice larger than $\bee$'s, where $\bT \bee$ is some coarse representation of $\bee$, say,
\be
    \lp \bT \bee \rp_U := \frac{1}{|\cT_U|} \sum_{u \in \cT_U} e_u\,,\qquad U \in \cN^c\,.
    \label{t}
\ee
$\bT$ is called the {\it aggregate type} operator \cite[\S 2]{cr_james}. (\ref{galerkin}) can be rewritten as
\be
   \mymin{\by^c} \left\{ \frac12 (\by^c)^{T} \bA^c \by^c - \by^c \bP^{T} \br \right\}\,,\qquad \br := \bA \bee\,.
   \label{min_ec_r}
\ee
For an ideal interpolation $\bP$ that satisfies $\bP \bT \bee = \bee$, (\ref{min_ec_r}) is minimized by $\by^c = \bT \bee$. A caliber-1 interpolation $\bP$ still satisfies $\bP \bT \bee \approx \bee$, but the first term in (\ref{min_ec_r}) is multiplied by the {\it energy inflation factor}
\be
    q(\bee) := \frac{E^c\left(\bT \bee\right)}{E\left(\bee\right)} = \frac{E\left(\bP \bT \bee\right)}{E\left(\bee\right)},
    \qquad E^c\left(\bee^c\right) := \frac12 \left(\bee^c\right)^{T} \bA^c \bee^c\,.
    \label{q}
\ee
Now (\ref{min_ec_r}) is minimized by $\by^c \approx q^{-1}\bT \bee$. As $\bee$ is not significantly changed by relaxation, its two-level Asymptotic Convergence Factor (ACF) will be $\rho \approx 1 - 1/q$. In Fig.~\ref{path}a, $q \approx 2$ and $\rho \approx .5$.

Several inflation remedies can be pursued:
\bi
    \item[(A)] Increase the $\bP$'s caliber and accuracy. This leads to fill-in troubles (\S \ref{interpolation_caliber}).
    \item[(B)] Accept an inferior two-level ACF of $1-1/q$ and increase the cycle index $\gamma$ to maintain it in a multilevel cycle. Unfortunately, not only does this increase complexity, the examples in \S \ref{energy_correction} demonstrate that $q$ can be arbitrarily large. This ACF cannot be improved by additional smoothing steps either, because it is governed by smooth mode convergence.
    \item[(C)] Correct the coarse level operator $\bA^c$ to match the fine level operator's energy during the setup phase. This option is considered in \S \ref{energy_correction}.
    \item[(D)] Modify the coarse level correction $\bP \bee^c$ to match the fine level error $\bee$ during the solve phase. This option is pursued in \S \ref{adaptive}.
\ei

\subsubsection{Flat Energy Correction}
\label{energy_correction}
\label{mu}
In this scheme, (\ref{galerkin}) is modified to
\be
    \bP^{T} \bA \bP \, \bee^c = \mu \bP^{T}\left(\bb - \bA \tbx\right)\,.
    \label{galerkin2}
\ee
The key question is how to choose $\mu$. Motivated by Fig.~\ref{path}a and its two-dimensional analogue, Braess used $\mu=1.8$, but his V-cycle convergence for 2-D grid graphs was mesh-independent only if a fixed number of levels were used per cycle, and if AMG was used as a preconditioner. In fact, no {\it predetermined global} factor exists that fits all error {\it corrections} in scenarios such as Fig.~\ref{path}b, because the coarse-level solution depends on all local inflation ratios, which vary among graph nodes.

On the other hand, a {\it local energy} correction factor $\mu$ does exist. Indeed, the fine-level and coarse-level quadratic energies are separable to nodal energies:
\begin{subequations} \label{nodal}
\begin{eqnarray}
    E(\bx) &=& \sum_{u \in \cN} E_u(\bx)\,,\quad\quad\,\,
    E_u(\bx) := -\frac12 \sum_{v: v \not = u} a_{uv} \left(x_u-x_v\right)^2\,, \label{nodal_fine} \\
    E^c(\bx^c) &=& \sum_{U \in \cN^c} E^c_U(\bx^c)\,,\quad
    E^c_U(\bx^c) := -\frac12 \sum_{V: V \not = U} a^c_{UV} \left(x^c_U-x^c_V\right)^2\,. \label{nodal_coarse}
\end{eqnarray}
\end{subequations}
Here $E_u(\bx)$ is the nodal energy at node $u$, and $E^c_U(\bx^c)$ is the nodal energy at aggregate $U$. The {\it local inflation factor} at aggregate $U$ is defined by
\be
    q_U(\bx) := \frac{E^c_U\left(\bT \bx\right)}{\sum_{u \in U} E_u\left(\bx\right)}\,.
    \label{q_local}
\ee
In principle, a {\it local} $\mu_U$ can be designed using our TVs to at least partially offset $q_U$; unfortunately, new difficulties arise (cf.~\S \ref{other_corrections}). Thus we chose to still scale the right-hand side by a global $\mu$, but {\it modify the aggregation so that $q_U(\bx) \lessapprox Q$ for all smooth vectors $\bx$ and all $U \in \cN^c$}, where $Q>1$ is a parameter. Under this condition a global factor is effective, whose optimal value minimizes the overall convergence factor:
\be
    \mu_{\text{opt}} = \argmin{\mu>0} \mymax{1 \leq q \leq Q} \left| 1 - \frac{\mu}{q} \right| =
        \argmin{\mu>0} \max \left\{ \left| 1 - \mu \right|, \left| 1 - \frac{\mu}{Q} \right| \right\} = \frac{2Q}{Q+1}\,.
    \label{optimal_mu}
\ee
In LAMG, we shoot for $Q=2$; hence $\mu_{\text{opt}}=\frac43$, and the expected ACF of smooth errors is $Q/(Q+1)=\frac13$.

\begin{figure}[htbp]
    \centering
    \begin{center}
    \begin{tabular}{cc}
      \includegraphics[height=1in]{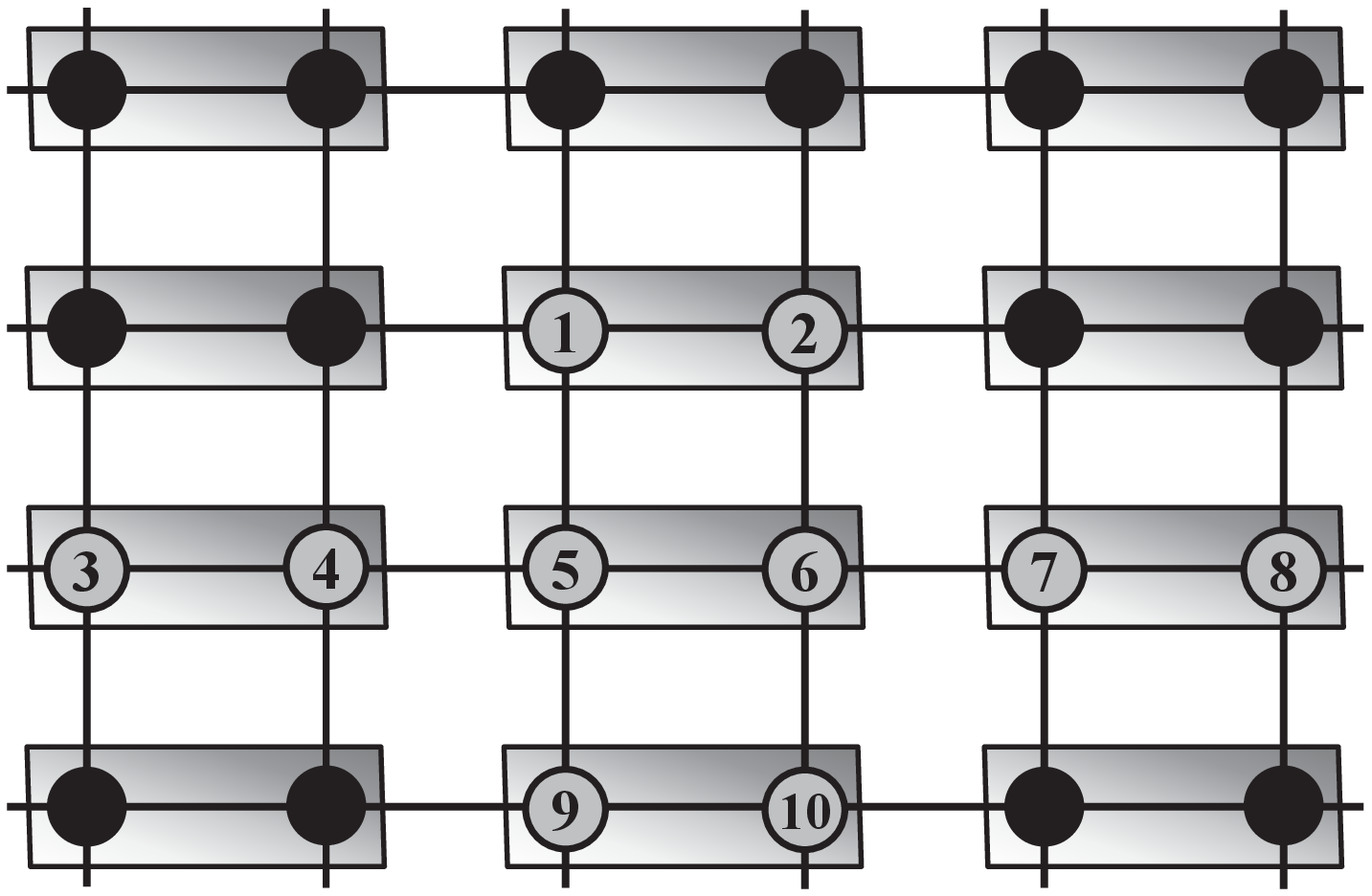}
      &
      \includegraphics[height=1in]{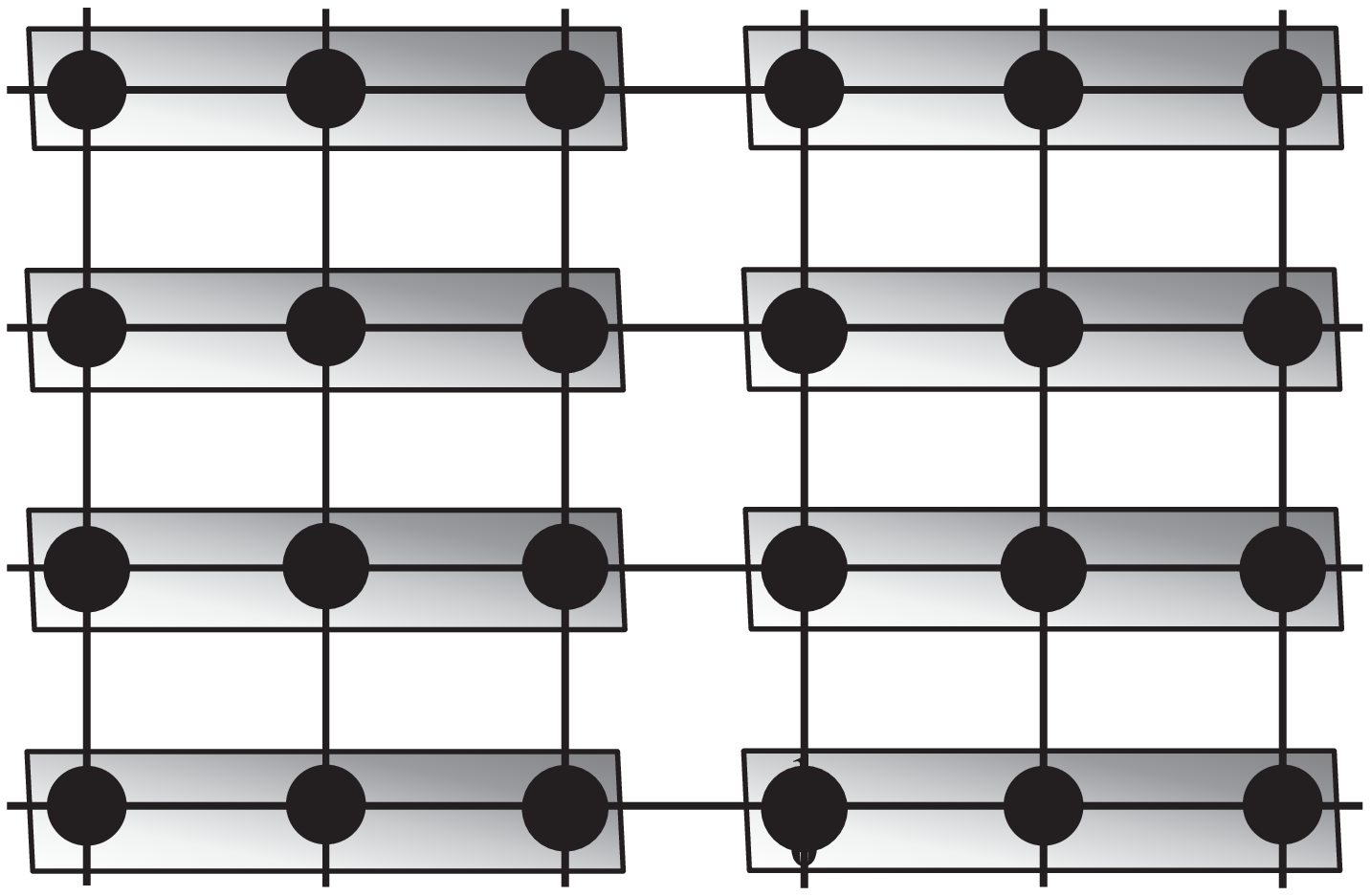}
      \\
      (a) $Q=2$ & (b) $Q=3$ \\
      \\
      \includegraphics[height=1in]{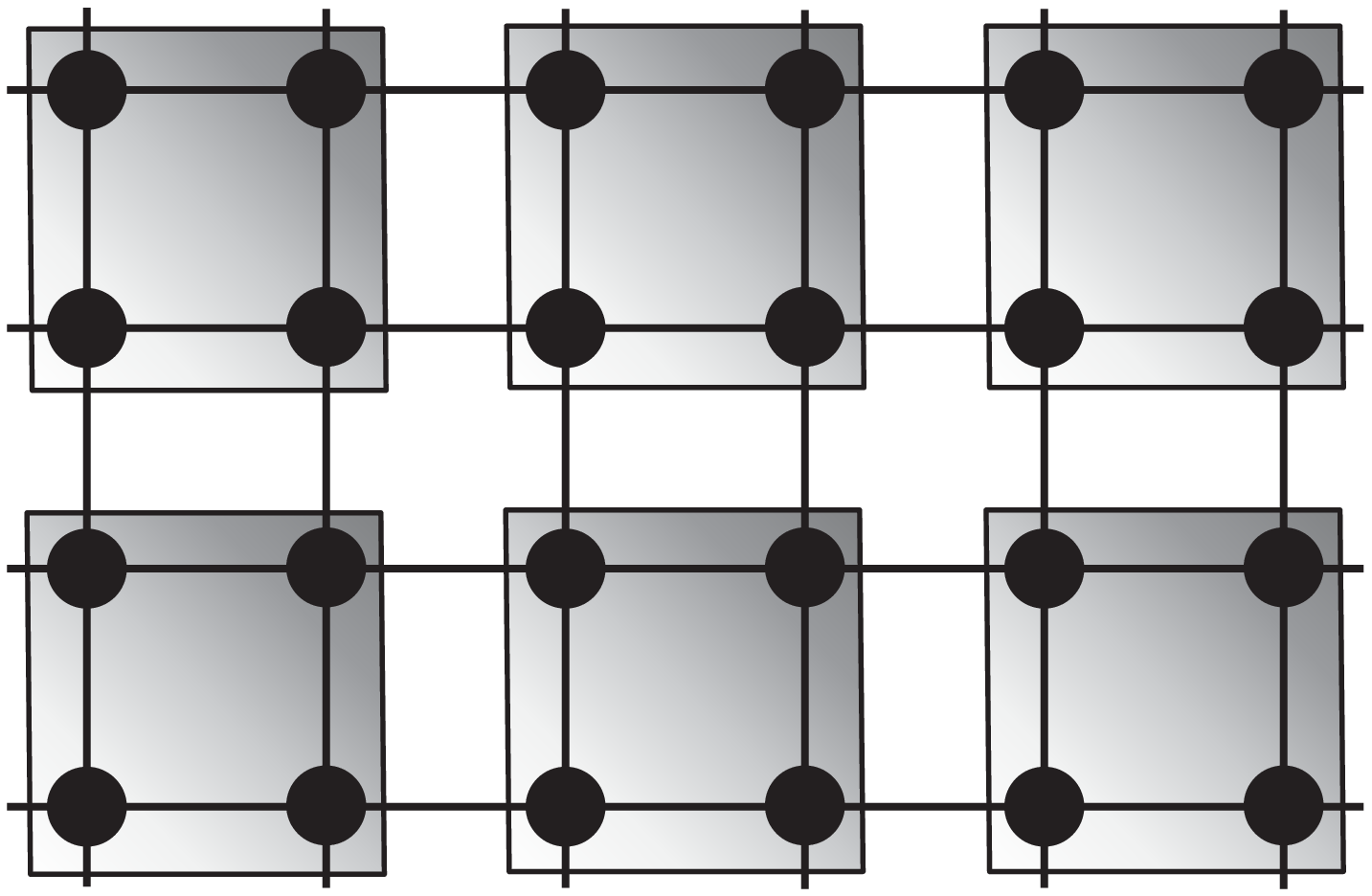}
      &
      \includegraphics[height=1in]{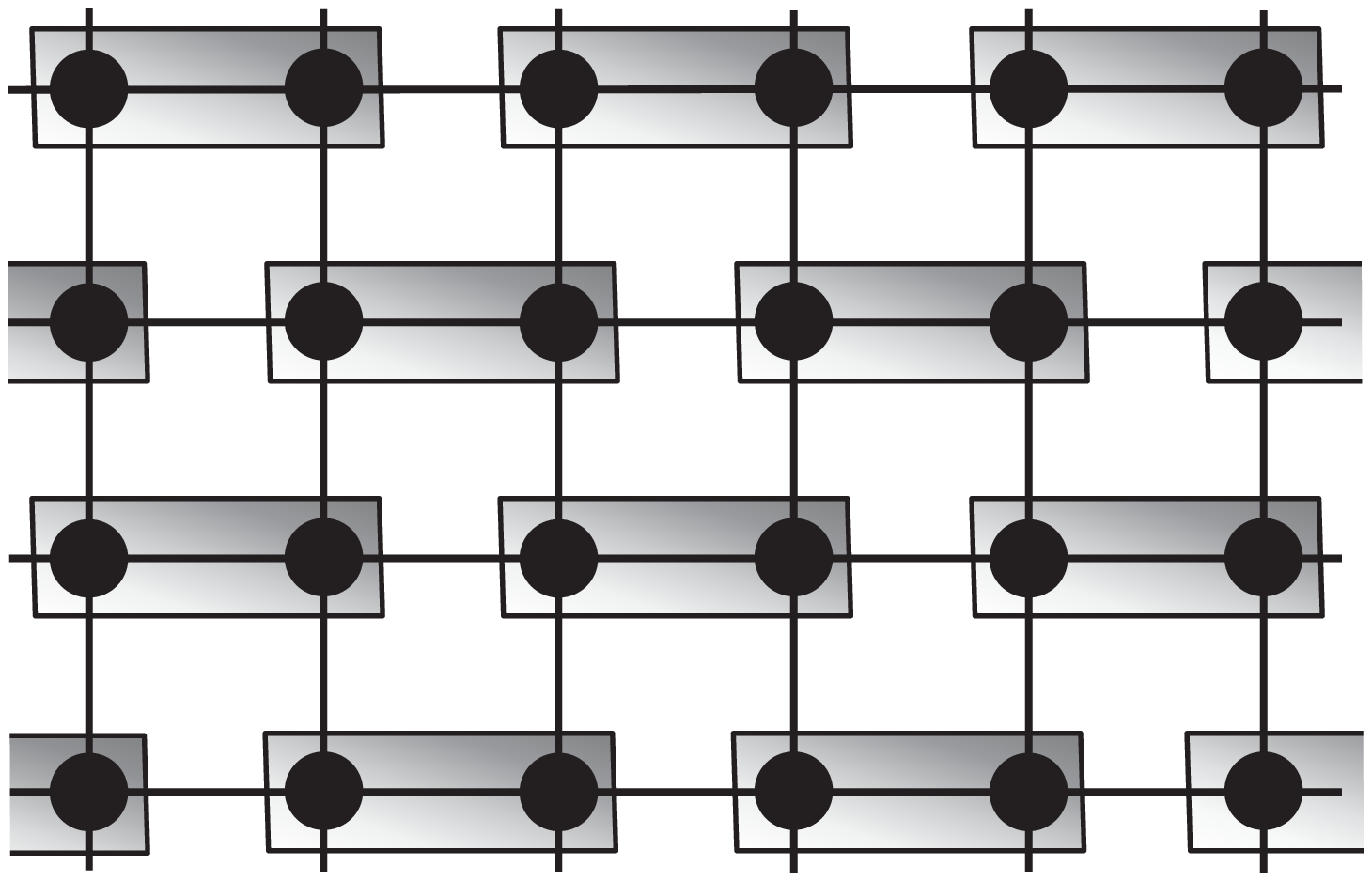}
      \\
      (c) $Q=2$ & (d) $Q=3$ \\
    \end{tabular}
    \end{center}
    \caption{Coarsening patterns. (a) 1:2 semi-coarsening. The energy ratio of aggregate $\{5,6\}$ depends on $\{x_u\}_{u=1}^{10}$. (b) 1:3 semi-coarsening.
    (c) 1:2 full coarsening. (d) Staggered semi-coarsening.}
    \label{grid2d}
\end{figure}

The worst energy ratio $Q := \max_{\bx} q_U(\bx)$ varies considerably with aggregate size, shape and alignment. Fig.~\ref{grid2d} depicts four constellations that may arise in an unweighted 2-D grid graph. (While these examples do not represent every scenario that can occur in graphs, they provide necessary conditions under which the algorithm must work.) Limiting the aggregate size to $2$, for instance, would not prevent case (d), whose $d$-dimensional analogue yields an unbounded $Q=d+1$. Fortuitously, we already possess the tool to signal and avoid bad aggregates: test vectors. The algorithm of \S \ref{agg_algorithm} is modified so that $u$ is only aggregated with a seed $s$ {\it if the local energy ratios of all test vectors are sufficiently small}.

Specifically, we compare the nodal energy $E_u$ before and after aggregation for each TV. Note that the nodal energy (\ref{nodal_fine}) is a quadratic in $x_u$ and $\{x_v\}_{v \in \cE_u}$. Define
$$
    E_u(\bx;y) := \frac12 a_{uu} y^2 - B_u(\bx) y + C_u(\bx)\,, B_u(\bx) := \sum_{v \in \cE_u} w_{uv} x_v\,, C_u(\bx) := \frac12 \sum_{v \in \cE_u} w_{uv} x_v^2\,,
$$
so $E_u(\bx) = E_u(\bx;x_u)$. The energy inflation that would occur upon aggregating $u$ with a seed $s$ is estimated by
\be
    q_{us} := \mymax{1\leq k \leq K} \frac{\mymin{y} E_u\left(\bx^{(k)};y\right)}{E_u\left(\bx^{(k)};x^{(k)}_s\right)}\,.
    \label{qut}
\ee
The numerator is the local energy after a temporary relaxation step is performed at $u$ (since the coarse-level correction is executed on a relaxed iterate during the cycle, this is the energy it aims to approximate; the papers \cite{bamg_markov, rbamg} use a similar idea). The denominator is the energy obtained when $x^{(k)}_u$ is set to $x^{(k)}_s$, simulating the caliber-1 aggregation; more accurate coarse-level energy estimates could be used, but we have not pursued them in the lean spirit of LAMG. We aggregate $u$ with the seed $s$ whose $c_{us}$ is minimal of all seeds $t$ with $q_{ut} \leq 2.5$; if none exist, $u$ is not aggregated at all. (Ratios slightly greater than the target $Q=2$ are accepted because TVs also contain high-energy modes, for which strict ratios are neither attainable nor necessary.) The complexity of the aggregation decision is $O(K |\cE_u|)$ \cite[\S 3.5.4]{lamg_arxiv}.

Low-degree elimination (\S \ref{elimination}) is advantageous because (a) it largely prevents worst case 1-D scenarios such as Fig.~\ref{path}b, where it is impossible to obtain low energy ratios without excessively increasing the coarsening ratio; (b) it increases the number of neighbors of $u$ and the chance of locating a seed $s$ with small energy inflation.

\subsubsection{Iterate Recombination}
\label{adaptive}

Instead of fixing $\mu$ by (\ref{galerkin2}), an effectively-{\it adaptive energy correction} is obtained by modifying the correction to smooth errors during the solution cycle. Let $l$ be any level such that $l+1$ is an \agg level. When the cycle switches from level $l-1$ to level $l$, $\vartheta$ sub-cycles are applied to $\bA^l \bx^l = \bb^l$, where $\vartheta$ is $1$ or $2$. We save the iterates $\bx^l_i$ obtained after the pre-relaxation of sub-cycle $i$, and, before switching back to level $l-1$, replace the final iterate $\bx^l$ by
\be
    \by^l = \bx^l + \balpha_1 \left( \bx^l_1 - \bx^l \right) + \cdots + \balpha_{\vartheta} \left( \bx^l_{\vartheta} - \bx^l \right)\,,
    \label{recomb}
\ee
where $\{\alpha_i\}_{i=1}^{\vartheta}$ are chosen so that $\|\bb^l - \bA^l \by^l\|_2$ is minimized (this is an $n_l \times \vartheta$ least-squares problem solved in $O(n_l)$ time). This {\it iterate recombination} \cite[\S 7.8.2]{trot} diminishes smooth errors $\bx^l_i-\bx^l$ that were not eliminated by $\kth{(l+1)}$-level corrections. Since the initial {\it residuals} obtained after interpolation from level $l+1$ are not smooth, the residual minimization is only effective after $\bx^l_i - \bx^l$ is smoothed. To maximize iterate smoothness, we perform more post- than pre-smoothing relaxations. The optimal splitting turned out to be a (1,2)-cycle; cf. \S \ref{algorithm}.

This acceleration is superior to CG because it is performed at all levels. Iterate recombination at coarse levels has been long recognized as an effective tool in the multigrid literature \cite[Remark~7.8.5]{trot}. In LAMG, recombination occurs more frequently at coarser levels because $\gamma>1$. Notay's K-cycle \cite{notay_etna, notay_nonsym} employs a similar multilevel CG acceleration, however with a much larger cycle index (up to $\gamma=4$), which increases the solver's complexity.

The aggregation is still modified here as in \S \ref{energy_correction}, to ensure small energy ratios and maximum reduction in the residual norm after recombination.

\section{The LAMG Algorithm}
\label{algorithm}

\subsection{Setup Phase}
\label{setup}
\begin{figure}[htbp]
    \centering
    \begin{center}
      \includegraphics[height=1.3in]{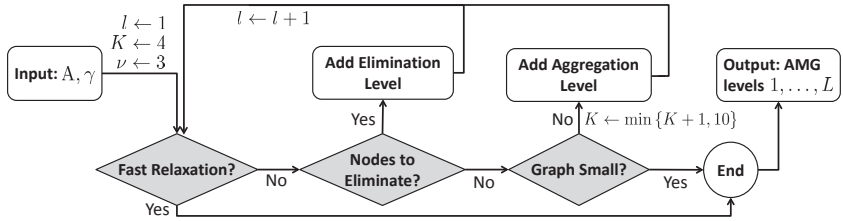}
    \end{center}
    \caption{LAMG setup phase flowchart.}
    \label{flowchart}
\end{figure}

The setup flow is depicted in Fig.~\ref{flowchart}. Its sole input is the cycle index $\gamma \geq 1$ to be employed at most levels of subsequent solution cycles. In our program, $\gamma=1.5$; this choice is discussed in \S \ref{solve}. The original problem ($l=1$) is repeatedly coarsened by either elimination or caliber-1 aggregation until the number of nodes drops below $150$, or until relaxation converges rapidly. We employ $K=4$ TVs at the finest level, and increase $K$ up to $10$ at coarser levels. Each TV is smoothed by $\nu=3$ relaxation sweeps.

\subsection{Solve Phase}
\label{solve}
The solve phase consists of multigrid cycles \cite[\S 1.4]{guide}. Each $l<L$ is assigned a cycle index $\gamma^l$ and pre- and post-relaxation sweep numbers $\nu^l_1,\nu^l_2$. If level $l+1$ is the result of elimination, $\gamma^l=1$ and $\nu^l_1=\nu^l_2=0$; otherwise,
\be
    \gamma^l :=
    \begin{cases}
        \gamma\,,& |\cE^l| > .1 |\cE|\,, \\
        \min\left\{2, .7\,|\cE^{l+1}|/|\cE^l| \right\}\,,& \text{otherwise} \,,
    \end{cases}\qquad
    \nu_1^l = 1\,,\quad\nu_2^l = 2\,.
    \label{cycle_index}
\ee
At fine levels, $\gamma=1.5$ is employed. This value is theoretically marginal for attaining a bounded multilevel ACF if the smoothest error two-level ACF is $\approx \frac13$ \cite[\S 6.2]{guide}, as implied by (\ref{optimal_mu}) for $Q=2$. Notwithstanding, worst-case energy ratios occur infrequently, and in practice a smaller ACF is obtained. This issue is further diminished by the adaptive energy correction. At coarse levels, $\gamma^l$ is increased to maximize error reduction while incurring a bounded work increase.

Three relaxation sweeps per level provide adequate smoothing, especially in light of the coarse-level correction's crudeness. The coarsest problem is solved by relaxation (if it is fast) or a direct solver on an augmented system \cite[\S 3.6.3]{lamg_arxiv}. Finally, (\ref{xcompat}) is enforced by subtracting the mean of $\bx$ from $\bx$ at the end of the cycle.

\begin{figure}[htbp]
    \centering
    \begin{center}
      \includegraphics[height=1.85in]{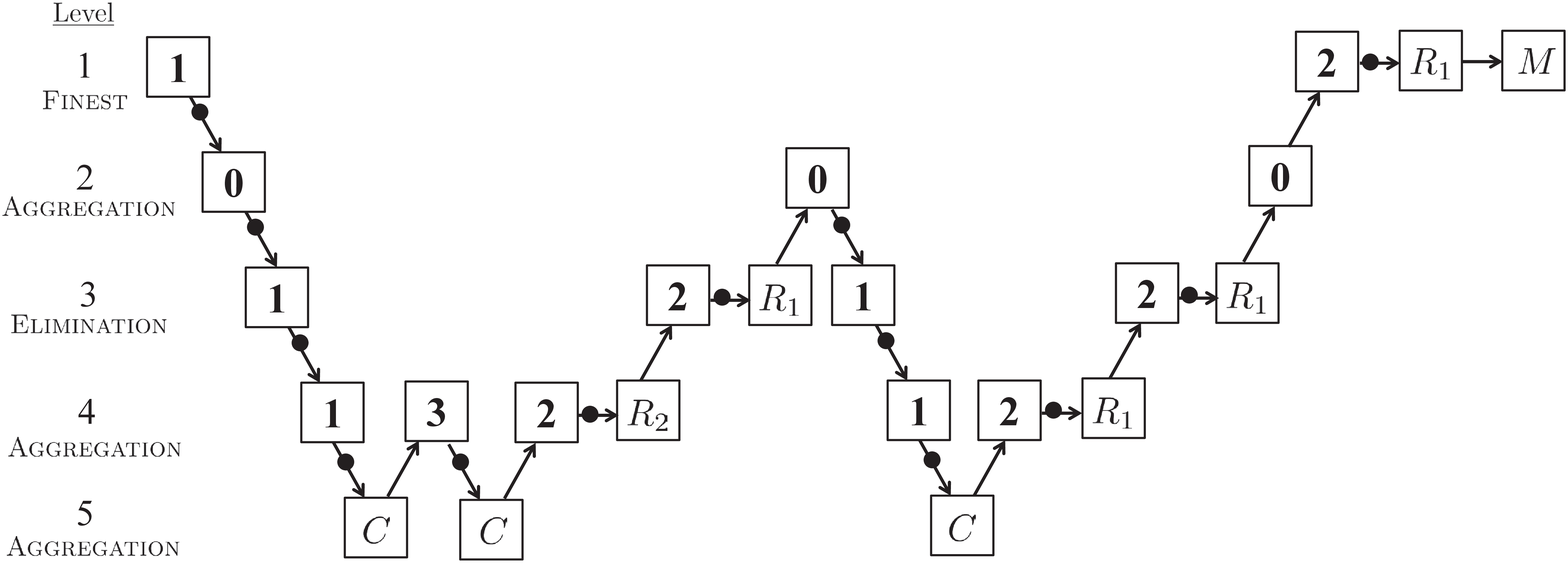}
    \end{center}
    \caption{A four-level cycle. A boxed number denotes a number of relaxations. $C$: the coarsest-level solver. $M$: subtracting the iterate mean. Down-arrows: right-hand side coarsening (\ref{elim}) or (\ref{galerkin2}). Up-arrows: coarse-level corrections \cite[Eq.~(3.4a)]{lamg_arxiv} or (\ref{correction}). $R_{\vartheta}$: a $(\vartheta+1)$-iterate recombination (\ref{recomb}). Iterates are saved at the black dots before coarsening.}
    \label{cycle}
\end{figure}

All cycle parameters are fixed: no fine tuning or parameter optimization is required for specific graphs. The total cycle work is equivalent to about $10$ relaxations.

\section{Numerical Results}
\label{results}
We provide supporting evidence for LAMG's practical efficiency for a wide range of graphs.

\subsection{Smorgasbord}
\label{smorgasbord}

An object-oriented \matlab 7.13 (R2011b) serial LAMG implementation was developed and is freely available online \cite{lamg_code}. The time-intensive functions were implemented in C and ported with the mex compiler \cite{matlab_davis}. It was tested on a diverse set of \numgraphs real-world graphs with up to \maxedges edges, collected from The University of Florida Sparse Matrix collection (UF) \cite{uf_collection}, C. Walshaw's graph partitioning archive \cite{walshaw}, I. Safro's MLogA results archive at Argonne National Laboratory \cite{mloga}, and the FTP site of the DIMACS Implementation Challenges \cite{dimacs}.

Graphs originated from a plethora of applications: airplane and car finite-element meshes; RF electrical circuits; combinatorial optimization; model reduction benchmarks; social networks; and web and biological networks. If the graph was directed, it was converted to undirected by summing the weights of both directions between each two nodes. Then, if it contained a large negative edge weight with $w_{uv} < -10^{-5} \sum_{v' \in \cE_u} |w_{uv'}|$, all weights were made positive by taking their absolute values. Finally, the Laplacian matrix $\bA$ was formed and used.

Runs were performed on Beagle, a 150 teraflops, 18,000-core Cray XE6 supercomputer at The University of Chicago (we only took advantage of parallelism by dividing the collection into equal parts, each of which ran a single AMD node with 2.2 GHz CPU and 32GB RAM). For each graph, a zero-sum random $\bb$ was generated. LAMG setup, followed by a linear solve that started from a random guess and proceeded until the residual $l_2$-norm was reduced by $10^{10}$. Six performance measures were computed:
\bi
    \item {\it Setup time per edge $\tsetup$}.
    \item {\it Solve time per edge per significant figure $\tsolve$}. If the residual norm after $i$ iterations was $r_i$ and $p$ iterations were executed, $\tsolve := t/(m \log_{10}(r_0/r_p))$, where $t$ was the solve time.
    \item {\it Total time per edge $\ttotal = \tsetup + 10 \tsolve$} to solve $\bA \bx=\bb$ to $10$ significant figures for a single $\bb$.
    \item Storage per edge.
    \item {\it Asymptotic convergence factor}, estimated by $(r_p/r_0)^{1/p}$.
    \item {\it Percentage spent on setup $\tsetup/\ttotal$}.
\ei

\begin{figure}[htbp]
    \centering
    \begin{center}
    \begin{tabular}{cc}
      \includegraphics[height=1.7in]{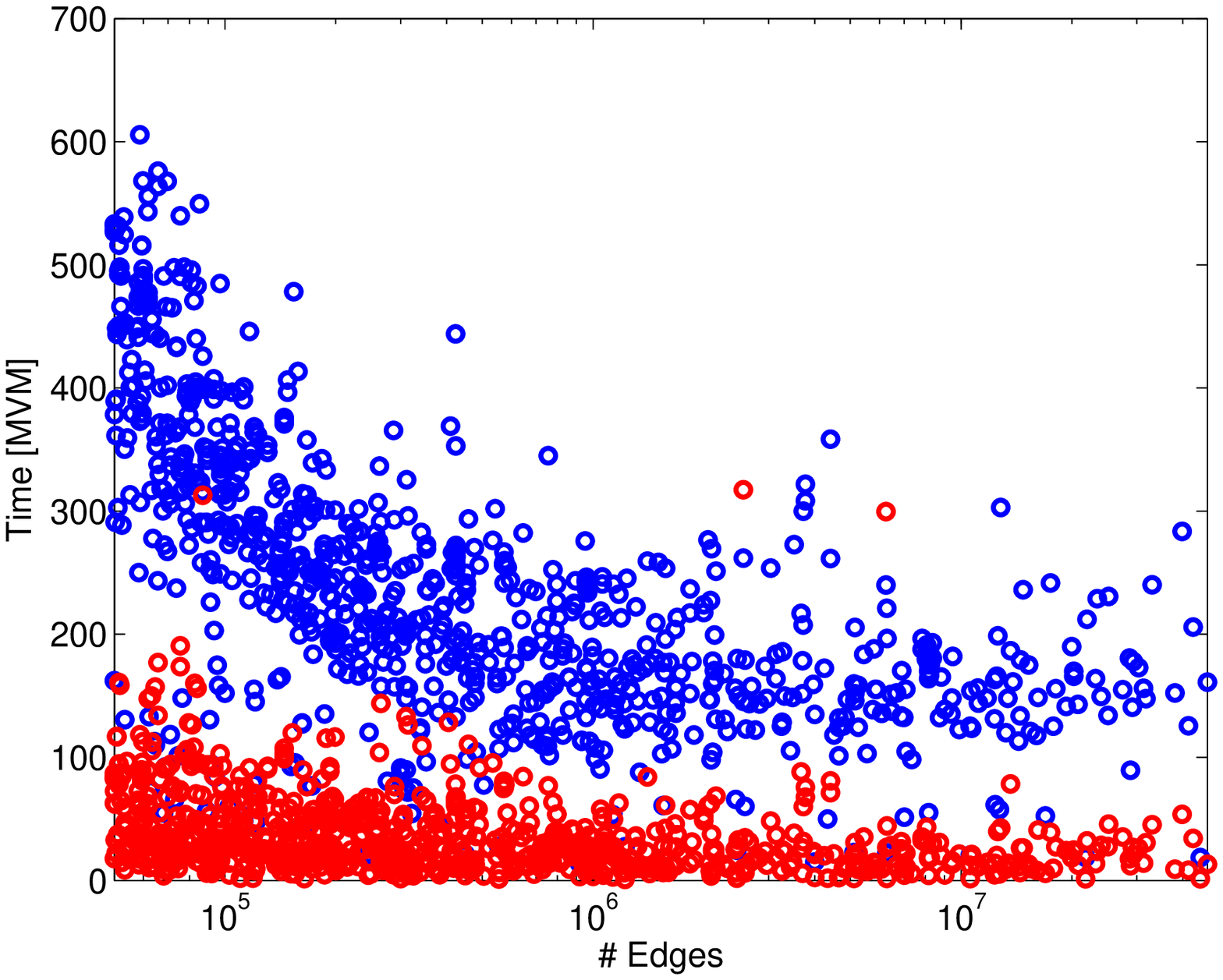}
      &
      \includegraphics[height=1.7in]{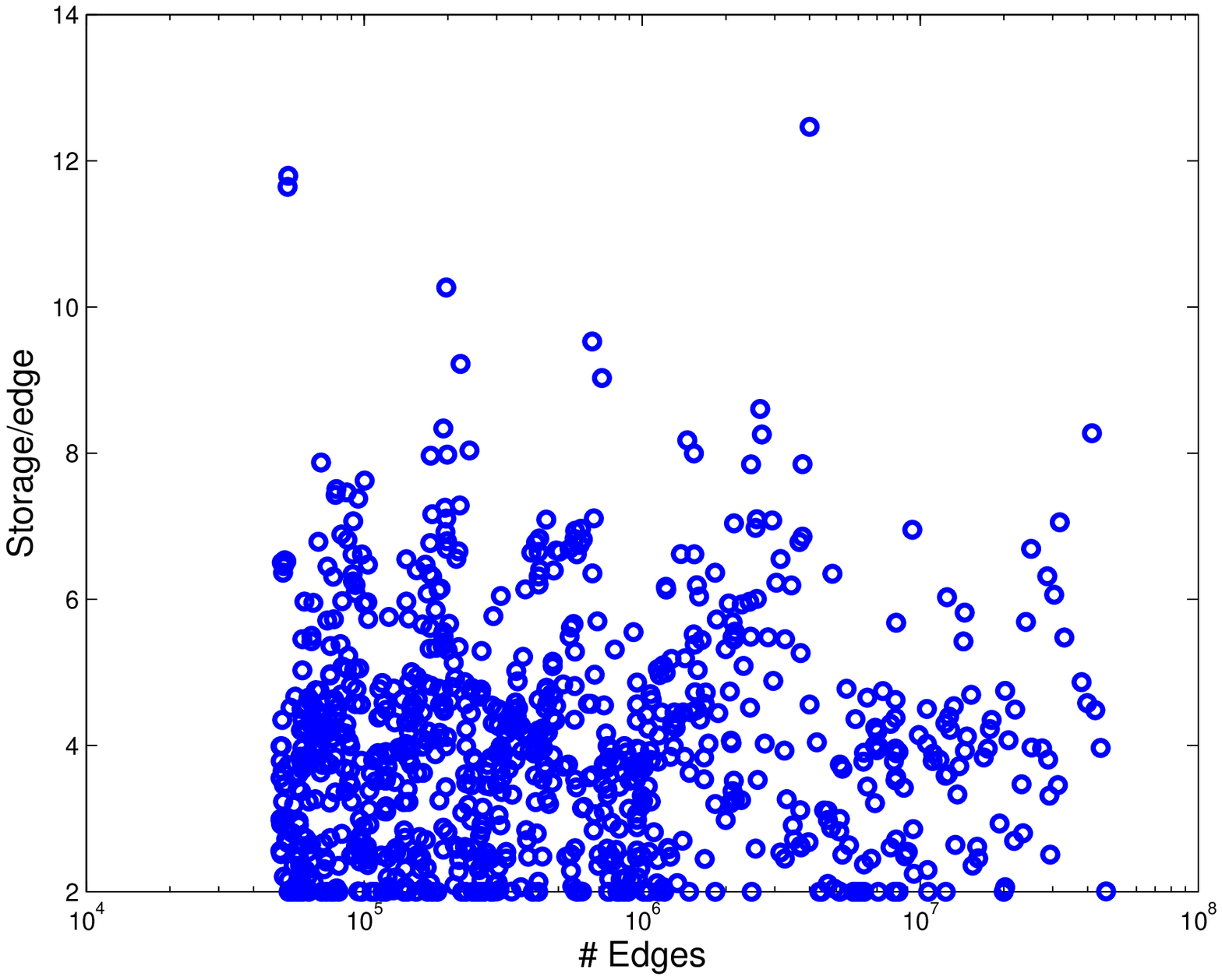}
      \\
      (a) & (b) \\
      \includegraphics[height=1.7in]{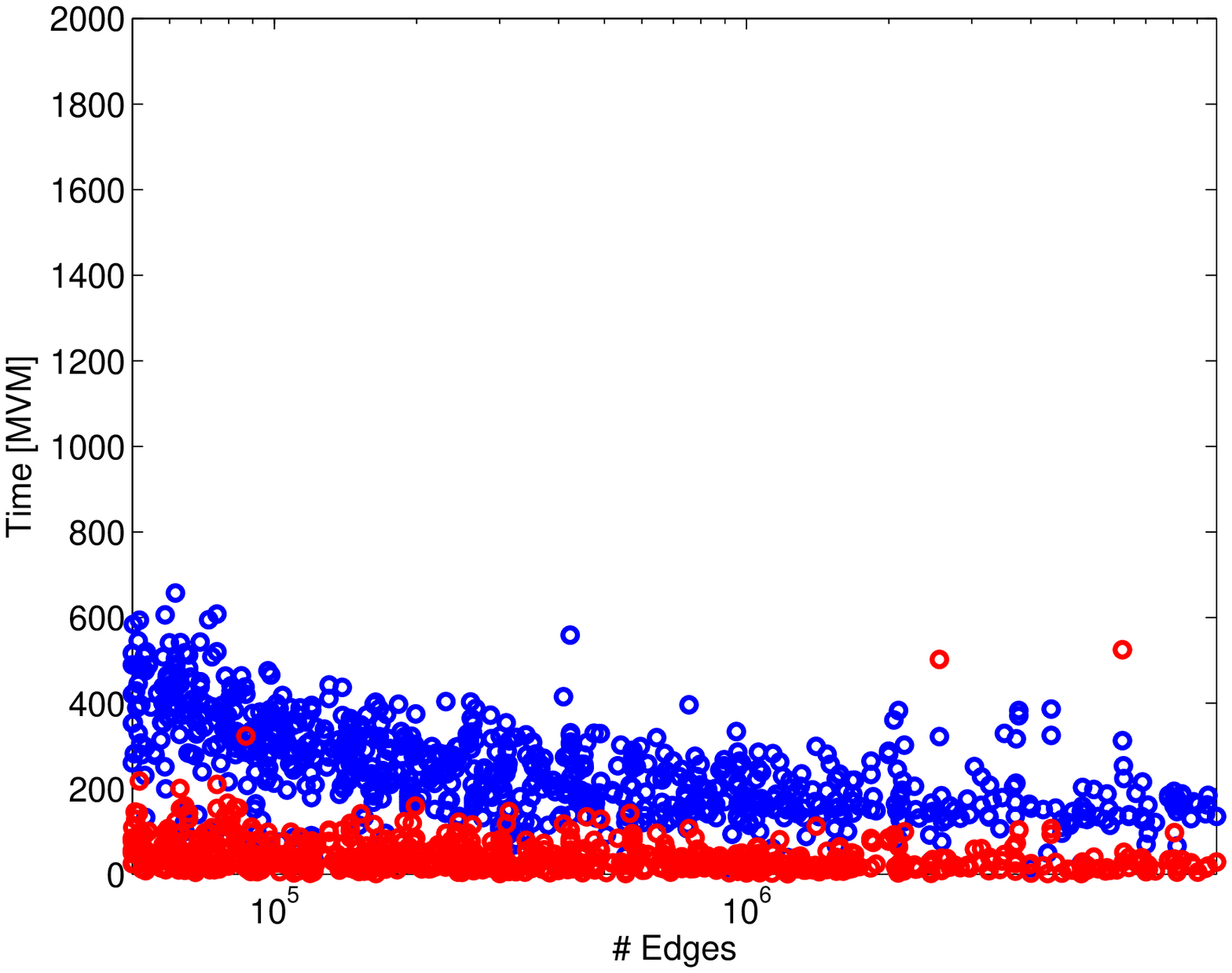}
      &
      \includegraphics[height=1.7in]{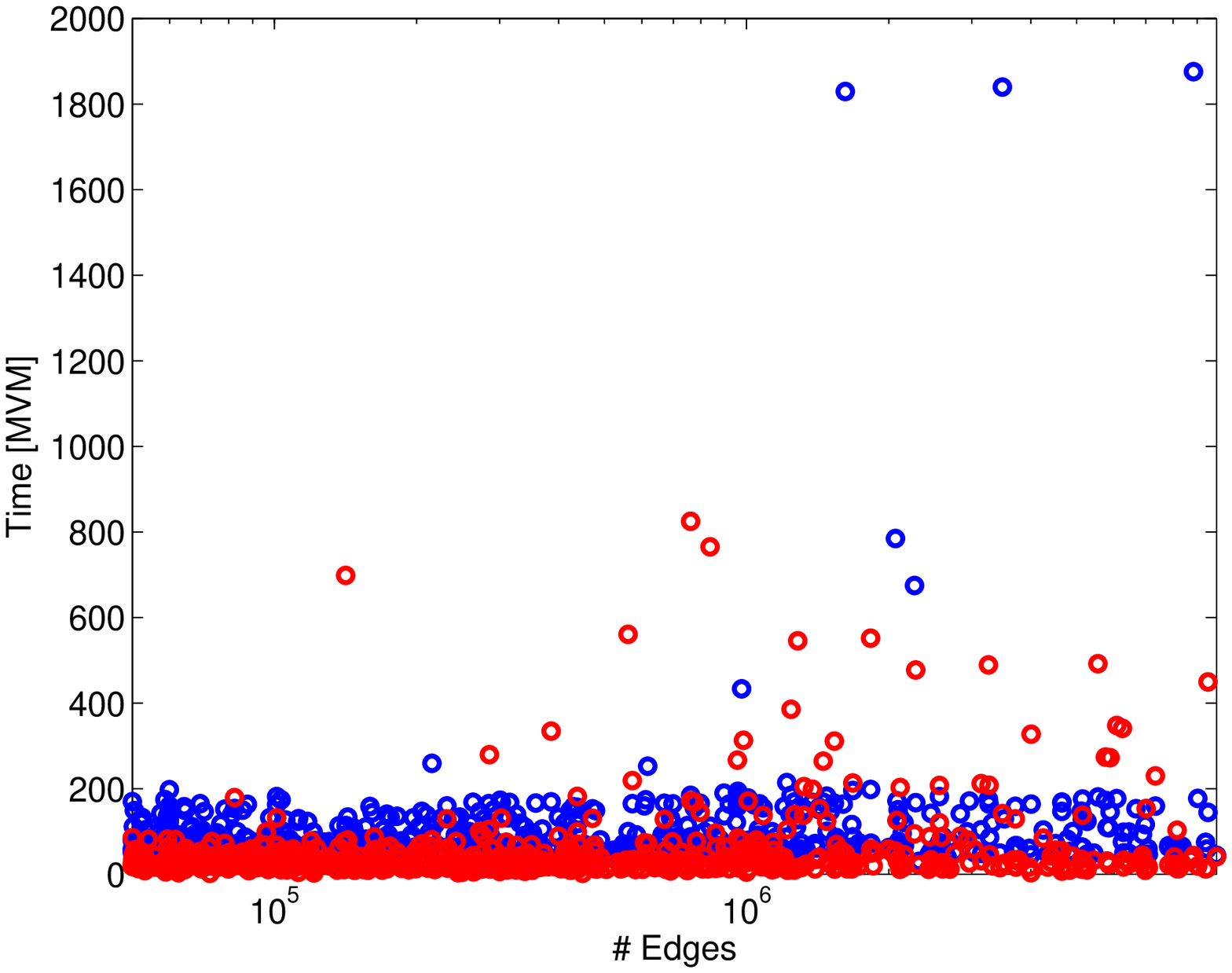}
      \\
      (c) & (d)
    \end{tabular}
    \end{center}
    \caption{(a) LAMG setup time (blue) and solve time (red) per edge on Beagle (up to $4.7 \times 10^7$ edges). (b) LAMG storage per edge on Beagle. (c) LAMG setup and solve time per edge on a Dell Inspiron (up to $10^7$ edges). (d) CMG setup and solve time per edge on a Dell Inspiron.}
    \label{times}
\end{figure}

LAMG scaled linearly with graph size: both $\tsetup$ and $\tsolve$ were approximately constant (Fig.~\ref{times}a; Table~\ref{times_avg}). Times were measured in terms of the most basic sparse matrix operation: a matrix-vector multiplication (MVM), because even MVM time scaled slightly superlinearly for $m \geq 5 \times 10^6$ due to loss of memory locality in the MATLAB compressed-column format \footnote{T. Davis, private communication.}. In wall clock time, the total time per edge was $5.6 \times 10^{-6}$ on average, i.e., LAMG performed a linear solve to $10$ significant figures at $178,000$ edges per second. The LAMG hierarchy required the equivalent of storing $\approx 4 m$ edges in memory (Fig.~\ref{times}b). Adaptive energy correction provided a $20\%$ speed up over a flat $\mu=\frac43$ and was thus employed in all reported experiments. The ACF was better than the expected $.33$ for flat correction (\S \ref{energy_correction}).

We compared LAMG with MATLAB's direct solver (the '$\backslash$' operator). Since the direct solver ran out of memory for many graphs with over $10^5$ edges, we did not include it in the plots. LAMG aims at robustness for a wide variety of graphs, and should be compared against solvers that do not often break down or require tuning, even if they are faster for a subset of the graphs (in analogy, many graphs could be solved much faster with a tailored geometric multigrid or classical AMG algorithm).

An advantage of iterative solvers over direct is their tunable solution accuracy $\ep$. Since $\bA$'s entries often incur measurement or modeling errors in applications, it does not make sense to solve $\bA\bx=\bb$ to more than 2--3 significant figures. Furthermore, a {\it single} multigrid cycle is typically sufficient to solve a nonlinear problems as well as time-dependent problems to the level of discretization errors \cite[Chaps.~7,15]{guide}.

We also compared LAMG against a \matlab implementation of CMG, a hybrid graph-theoretic-AMG solver \cite{cmg}. Since CMG doesn't run yet on the Cray architecture, experiments were performed on a smaller 64-bit Dell Inspiron 580 (3.2 GHz CPU; 8GB RAM) for $2668$ graphs with up to $10^7$ edges. Both algorithms successfully solved all graphs. Solve times were similar, while LAMG's setup time was thrice larger than CMG's. On the other hand, LAMG was much more robust than CMG: it had only 3 outliers whose solve time was large, as opposed to 26 CMG outliers whose relative magnitude was much larger (4 in setup and 22 in solve; see Fig.~\ref{times}c-d and Table~\ref{outliers}).

\begin{table}[htbp]
\centering \footnotesize
\begin{tabular}{|c||c|c||c|c|c|c|}
\hline
\multirow{2}{*}{Measure} &
\multicolumn{2}{c||}{LAMG (Beagle)} &
\multicolumn{2}{c|}{LAMG (Dell)}   &
\multicolumn{2}{c|}{CMG (Dell)}    \\ \cline{2-7}
& Median & Mean $\pm$ Std. & Median & Mean $\pm$ Std. & Median & Mean $\pm$ Std. \\ \hline \hline
$\ttotal$  & $482.9$ & $585.4 \pm 406 $ & $558.1$ & $680.5 \pm 497 $ & $342.2$ & $582.6 \pm 1088$ \\ \hline
$\tsetup$  & $199.6$ & $222.5 \pm 108 $ & $234.9$ & $248.4 \pm 113 $ & $66.3$  & $89.3  \pm 111  $ \\ \hline
$\tsolve$  & $27.3$  & $36.3  \pm 33  $ & $31.6$  & $43.2  \pm 42.3$ & $25.5$  & $47.9  \pm 106  $ \\ \hline
ACF        & $.107$  & $.128  \pm .12 $ & $.112$  & $.132  \pm .11 $ & $.500$  & $.495  \pm .21 $ \\ \hline
\%Setup    & $43.8\%$  & $43.7\%  \pm 13\%  $ & $42.4\%$  & $43.4\%  \pm 13\%  $ & $21.5\%$  & $22.9\%  \pm 12\%$ \\ \hline
\end{tabular}
\caption{Left column: median and mean LAMG performance on the Beagle Cray for $892$ graphs with $50,000 \leq m \leq 4.7 \times 10^7$. Middle and right: LAMG vs. CMG performance on a Dell Inspiron for $794$ graphs with $50,000 \leq m \leq 10^7$. Times are measured in matrix-vector multiplications.}
\label{times_avg}
\end{table}

\begin{table}[htbp]
\footnotesize
\centering
\begin{tabular}{|p{1.3in}|c|c||c|c|c||c|c|}
\hline
Name & $n$ & $m$ &
\multicolumn{3}{c|}{LAMG (Dell)}   &
\multicolumn{2}{c|}{CMG (Dell)}    \\ \cline{4-8}
                       &&& ACF       & $\tsetup$ & $\tsolve$ & $\tsetup$ & $\tsolve$ \\ \hline \hline
Ill-conditioned Stokes   & $  20896$ & $  87010$ & $\mathbf{.66}$  & $ 420$ & $ \mathbf{323}$ & $  66$ & $  25$ \\ \hline
Large basis              & $ 440020$ & $2560040$ & $\mathbf{.88}$  & $ 322$ & $ \mathbf{502}$ & $  70$ & $ 107$ \\ \hline
RF circuit simulation    & $4690002$ & $6251251$ & $\mathbf{.72}$  & $ 312$ & $ \mathbf{525}$ & $  65$ & $ 304$ \\ \hline \hline
Law citation network     & $ 925340$ & $6675561$ & $.24$  & $ 169$ & $  15$ & $ 152$ & $\mathbf{2037}$ \\ \hline
Berkeley-Stanford web    & $ 512501$ & $3480880$ & $.17$  & $ 168$ & $  18$ & $\mathbf{1585}$ & $ 126$ \\ \hline
Molecule pseudopotential & $ 268096$ & $8833823$ & $.13$  & $ 167$ & $  21$ & $\mathbf{1879}$ & $  41$ \\ \hline
\end{tabular}
\caption{Top section: the three LAMG outliers. Bottom section: three of CMG's 26 outliers. Times are measured in matrix-vector multiplications.}
\label{outliers}
\end{table}

\begin{figure}[htbp]
    \centering
    \begin{center}
    \begin{tabular}{cc}
      \includegraphics[height=1.35in]{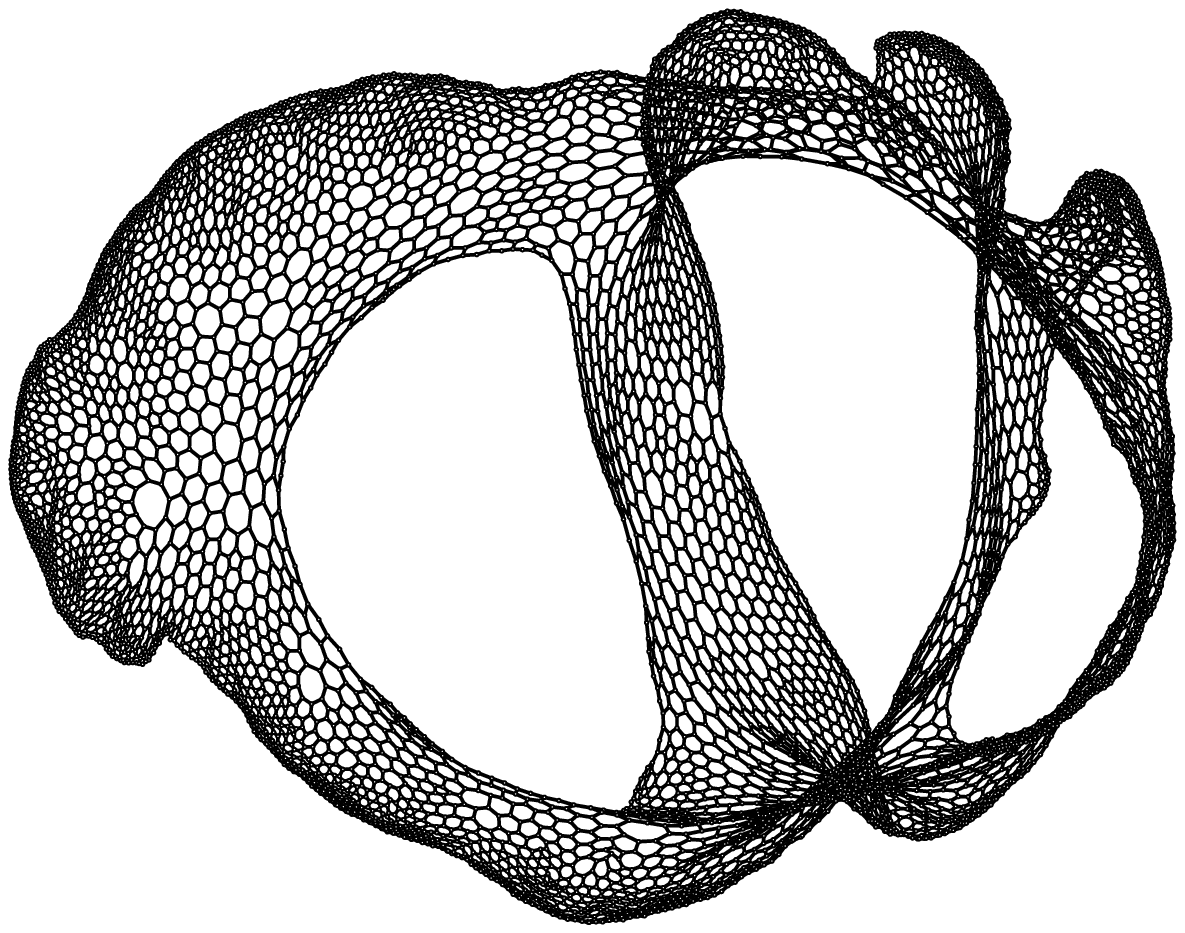}
      &
      \includegraphics[height=1.35in]{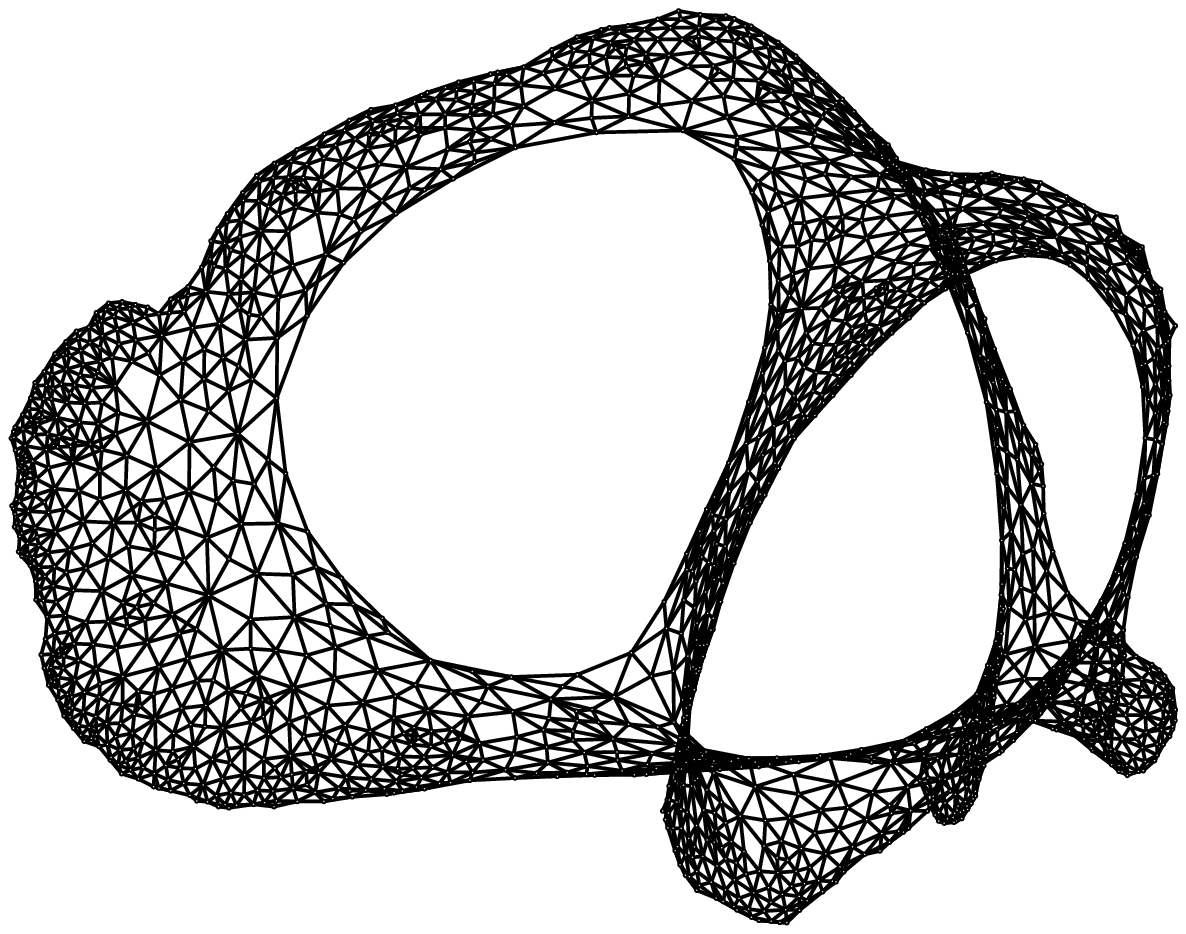}
      \\
      $G^1$ & $G^3$ \\
      \includegraphics[height=1.35in]{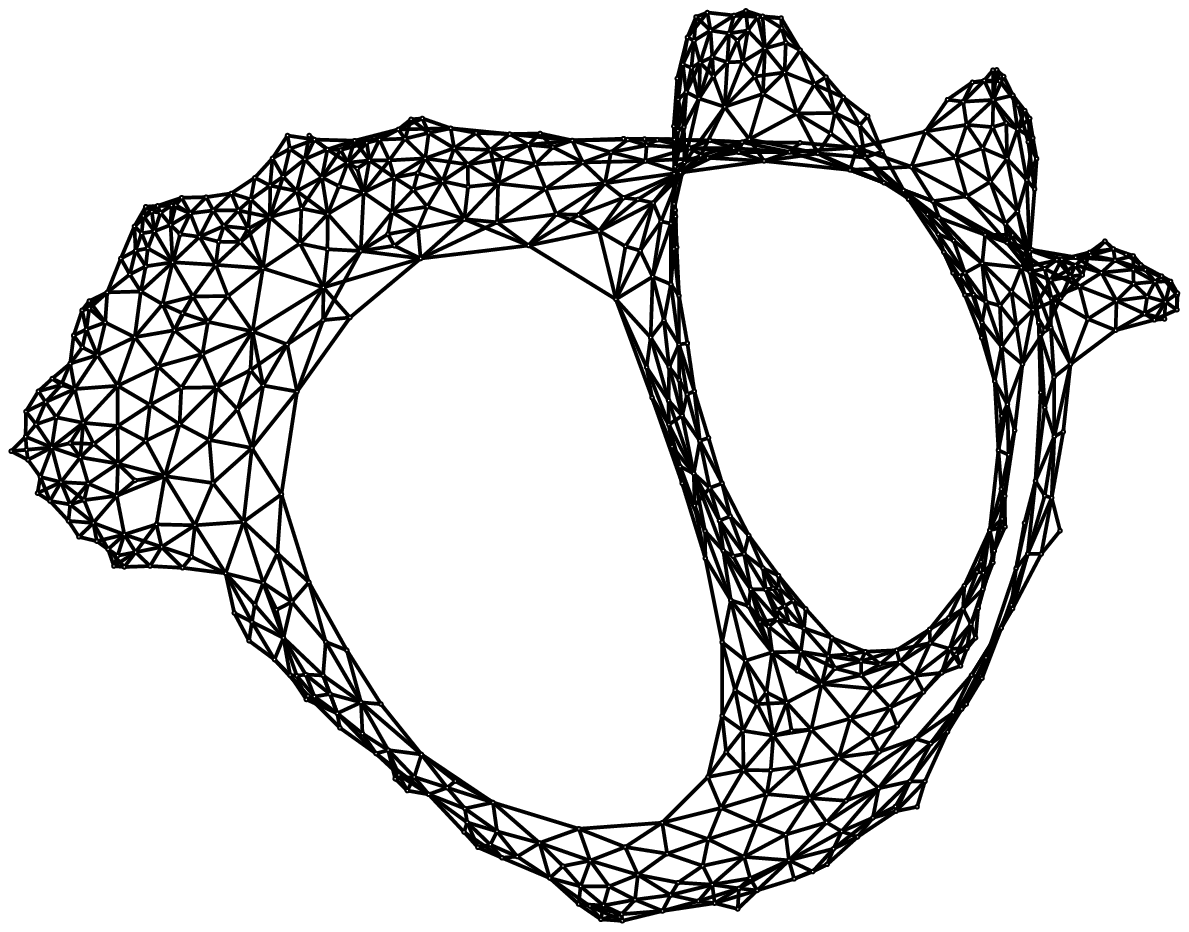}
      &
      \includegraphics[height=1.35in]{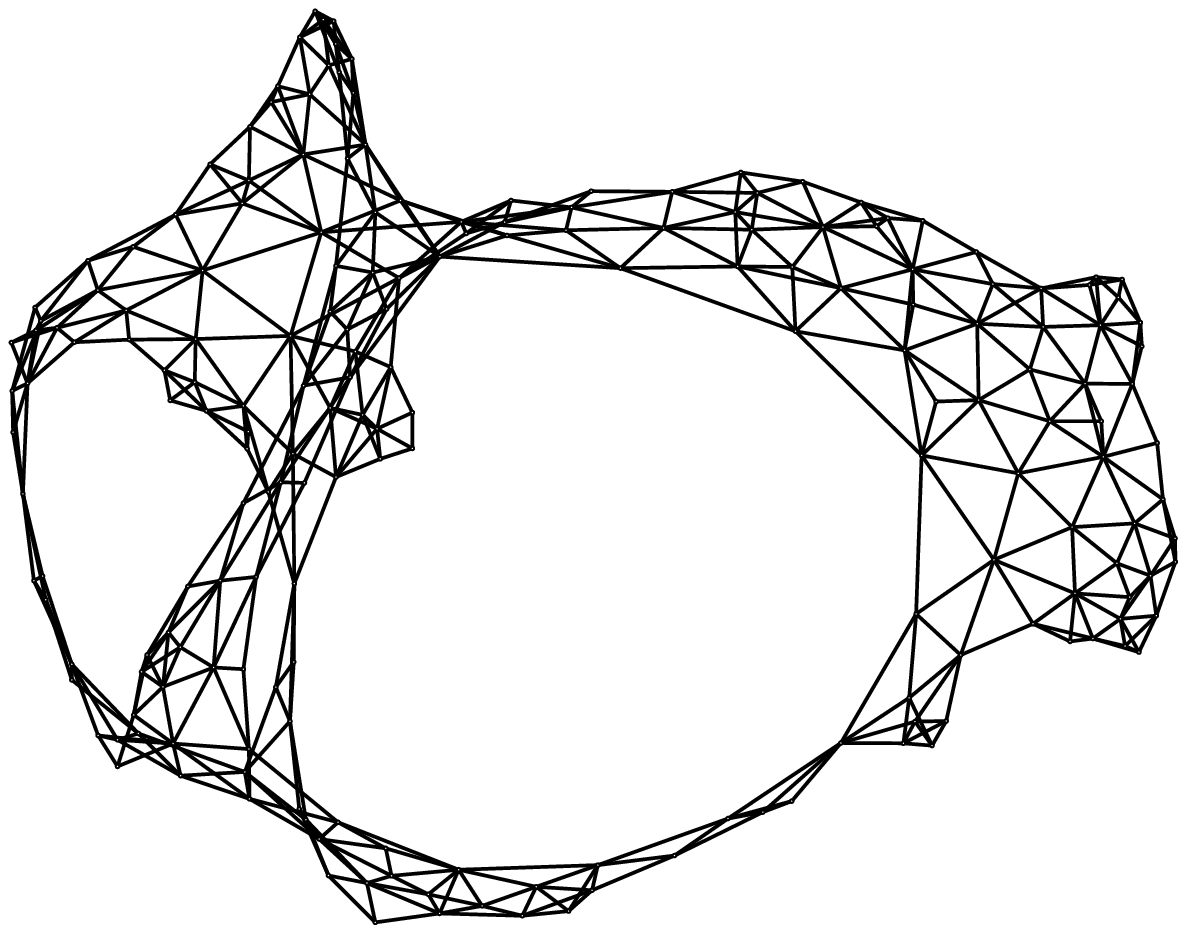}
      \\
      $G^5$ & $G^7$
    \end{tabular}
    \end{center}
    \caption{The four finest aggregation levels for the UF 2-D airfoil finite-element planar graph
    {\tt AG-Monien/airfoil1-dual}. Graphs were drawn using GraphViz with the SFDP algorithm \cite{graphviz}.}
    \label{airfoil_graphs}
\end{figure}
\begin{figure}[htbp]
    \centering
    \begin{center}
    \begin{tabular}{cc}
      \includegraphics[width=1.7in,height=1.35in]{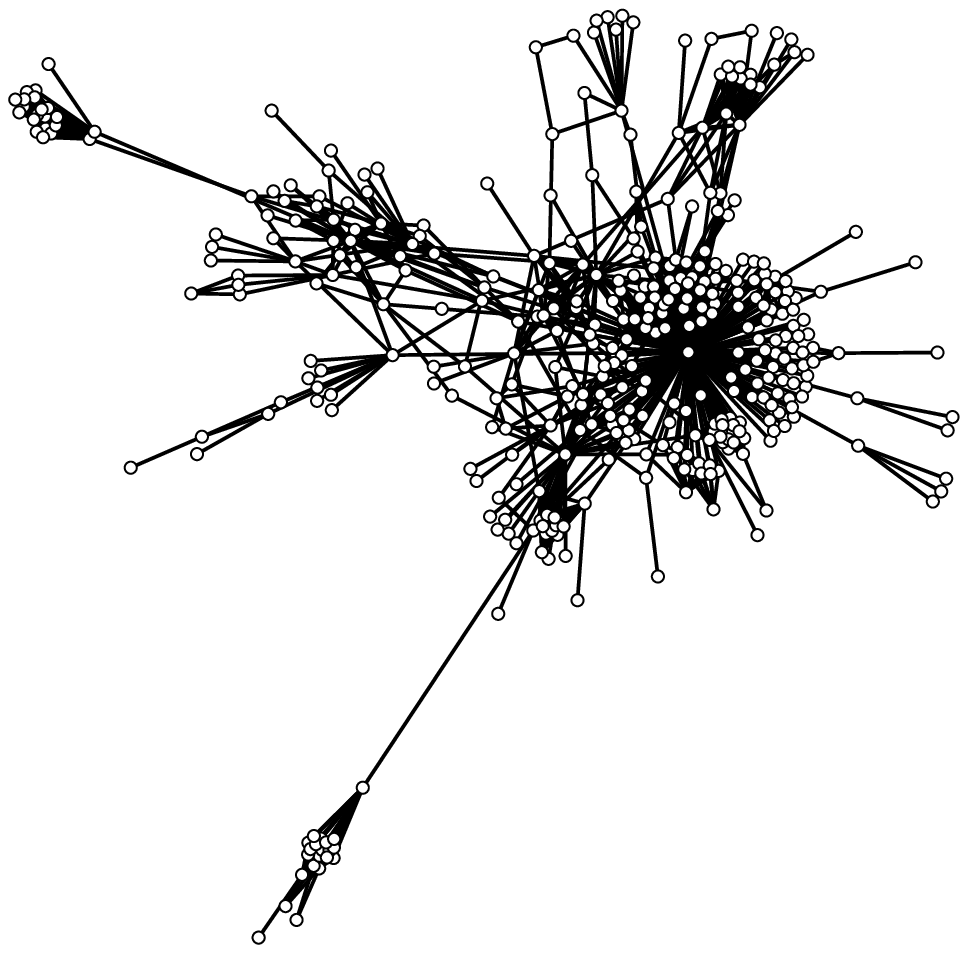}
      &
      \includegraphics[width=1.7in,height=1.35in]{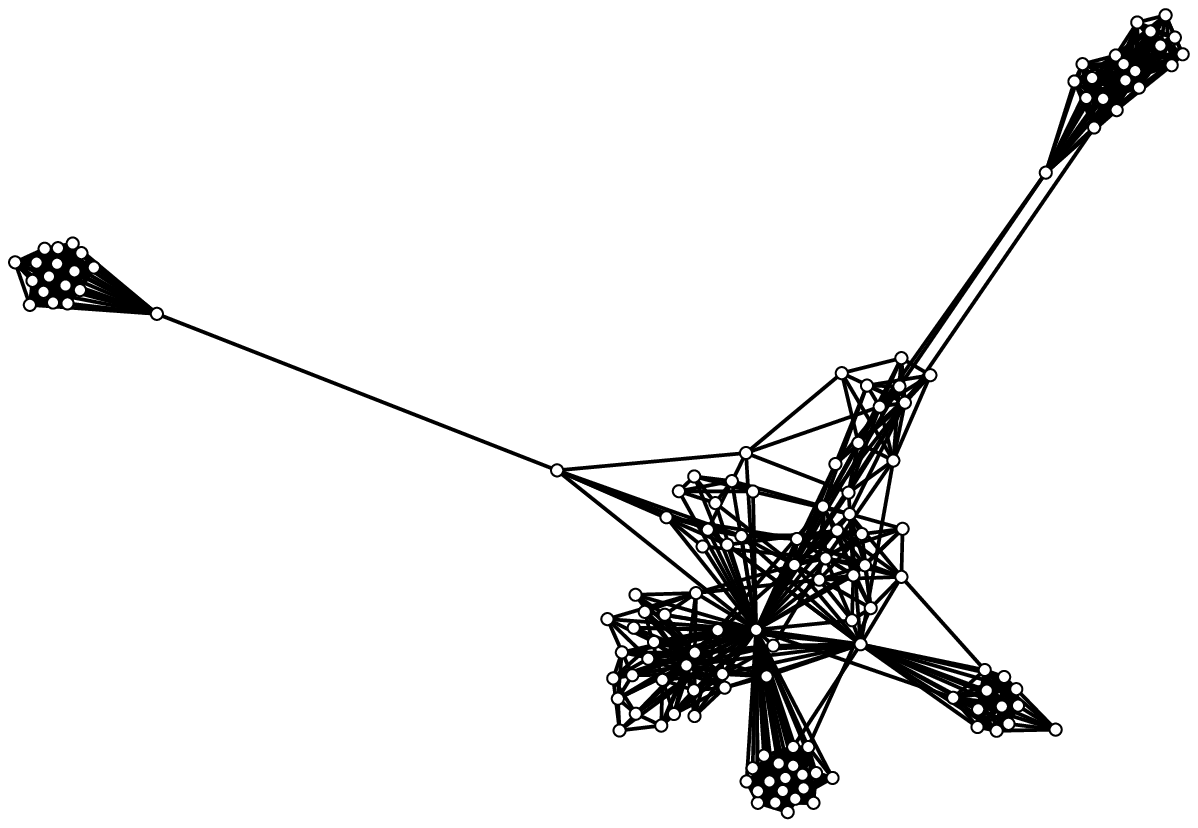}
      \\
      $G^1$ & $G^2$ \\
      \includegraphics[width=1.7in,height=1.35in]{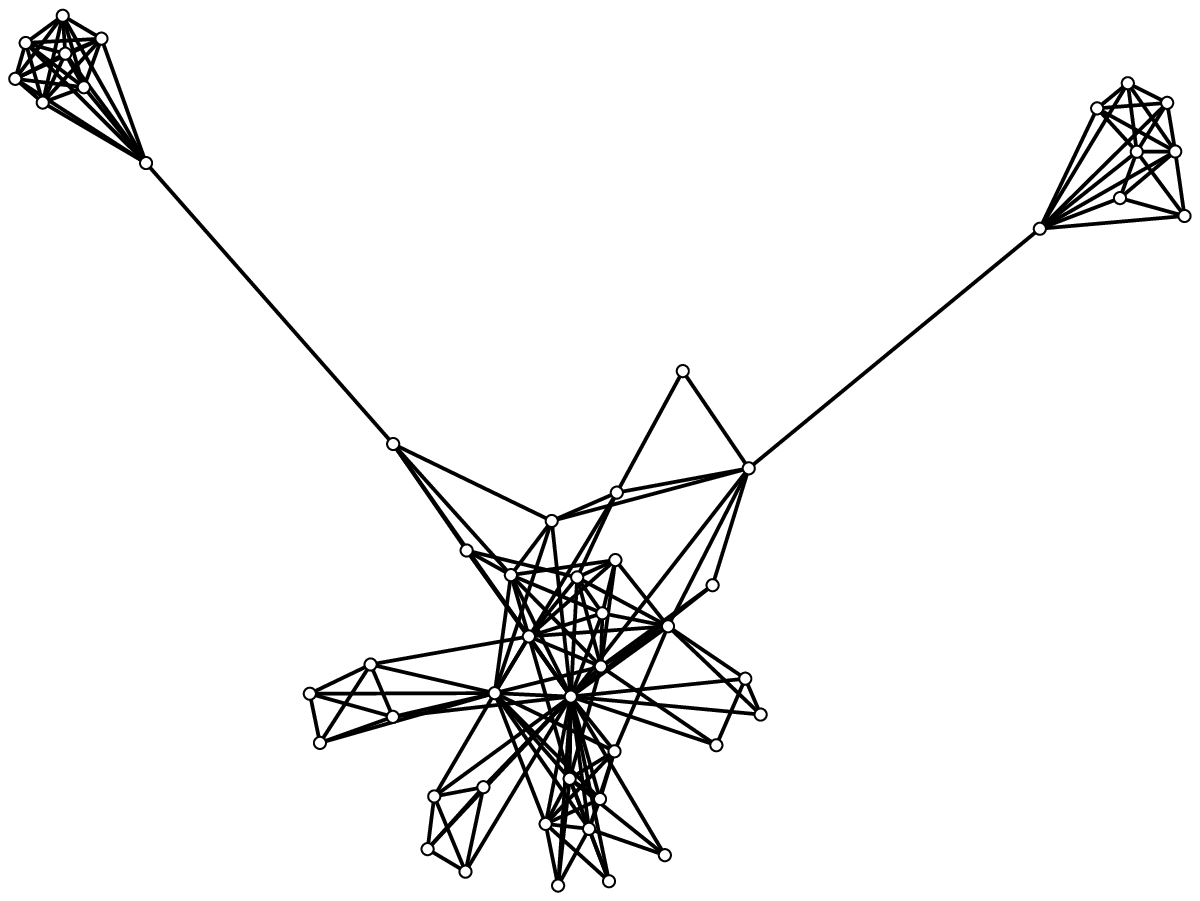}
      &
      \includegraphics[width=1.7in,height=1.2in]{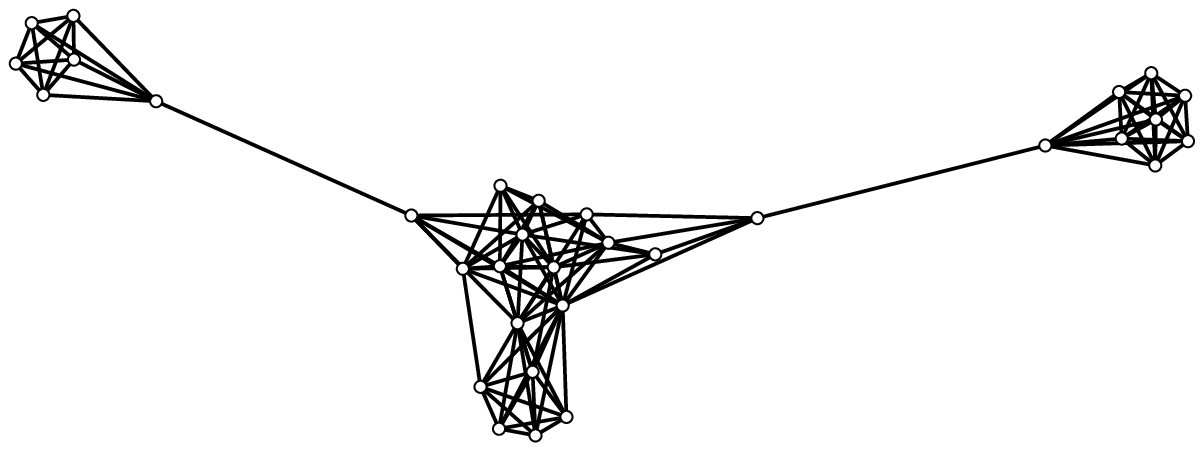}
      \\
      $G^3$ & $G^4$
    \end{tabular}
    \end{center}
    \caption{The four finest levels for the UF Harvard 500 non-planar web graph \cite{matlab_moler}.}
    \label{harvard500_graphs}
\end{figure}

\subsubsection{LAMG's Outliers}
The three solve-time outliers (Table~\ref{outliers}) were characterized by a large portion of small edge weights, which were carried over to all coarse matrices and increased coarsening ratios. Since LAMG's work was controlled by decreasing $\gamma$ via (\ref{cycle_index}), a slower cycle resulted. We plan to improve those cases in the future by appropriately ignoring weak edges at each level; cf.~\S \ref{improvements}.

\subsection{Grids with Negative Weights}
\label{negative_weights}

\begin{figure}[htbp]
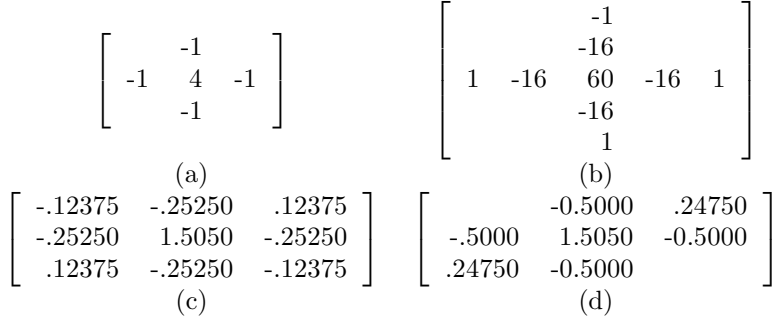

    \centering
    \begin{center}
    \begin{tabular}{cc}
   $\left[
        \begin{tabular}{rrr}
                   &   -1  &      \\
               -1  &    4  &  -1    \\
                   &   -1  &
        \end{tabular}
    \right]$
    &
    $\left[
        \begin{tabular}{rrrrr}
            &       &   -1  &       & \\
            &       &  -16  &       & \\
        1   & -16   &   60  & -16   & 1 \\
            &       &  -16  &       & \\
            &       &    1  &       &
        \end{tabular}
    \right]$
      \\
      (a) & (b) \\
    $\left[
        \begin{tabular}{rrr}
        -.12375   &   -.25250  &   .12375  \\
        -.25250   &    1.5050  &  -.25250   \\
         .12375   &   -.25250  &  -.12375
        \end{tabular}
    \right]$
    &
    $\left[
        \begin{tabular}{rrr}
                   & -0.5000  &  .24750  \\
          -.5000   &  1.5050  & -0.5000  \\
          .24750   & -0.5000  &
        \end{tabular}
    \right]$
      \\
      (c) & (d) \\
    \end{tabular}
    \end{center}
    \caption{Stencils of grid Laplacians with negative weights. Entries are normalized to a meshsize-independent sum and rounded to five significant figures.}
    \label{negative_stencils}
\end{figure}

Unlike CMG, LAMG is not restricted to diagonally-dominant systems, and can also be applied to some graphs with negative edge weights $w_{uv}$, as long as the Laplacian matrix is (or is very close to being) positive semi-definite. To demonstrate this capability, we tested LAMG on the following SPS 2-D grid Laplacians, whose stencils are depicted in Fig.~\ref{negative_stencils}:
\bi
    \item[(a)] The standard 5-point finite-difference discretization of $U_{xx} + U_{yy}$ on the unit square with Neumann boundary conditions.
    \item[(b)] The 13-point $\kth{4}$-order finite-difference stencil of $U_{xx} + U_{yy}$.
    \item[(c)] The discretized anisotropic-rotated Laplace operator
    \be
        \left(\cos^2\alpha + \ep \sin^2\alpha\right) U_{xx} + (1-\ep) \sin(2\alpha) U_{xy}
        + \left(\ep \cos^2\alpha + \sin^2\alpha\right) U_{yy}\,,
        \label{anis_rot}
    \ee
    with $\alpha=-\pi/4$, $\ep=10^{-2}$, standard 5-point stencil of $U_{xx}$,$U_{yy}$, and an alignment-agnostic cross-term
    $$
        U_{xy} \approx
           \frac{1}{4 h^2}
    \left[
        \begin{tabular}{rrr}
              -1   &    0  &   1   \\
               0   &    0  &   0    \\
               1   &    0  &  -1
        \end{tabular}
    \right]\,,$$
    where $h$ is the grid meshsize. Neumann boundary conditions were used.
    \item[(d)] The same as (c), but aligning $U_{xy}$ with the northeast and southwest neighbors:
    $$
        U_{xy} \approx
           \frac{1}{2 h^2}
    \left[
        \begin{tabular}{rrr}
               0   &   -1  &   1   \\
              -1   &    2  &  -1    \\
               1   &   -1  &   0
        \end{tabular}
        \right]\,.$$
\ei
Problems (c) and (d) are bad discretizations that do not align with the characteristic direction of (\ref{anis_rot}), and are considered hard for AMG \cite{alg_distance_anis}. Performance figures for the Dell Inspiron are given in Table~\ref{results_negative}.

\begin{table}[htbp]
\centering\footnotesize
\begin{tabular}{|p{3.3cm}|c|c|c|c|c|c|}
\hline
Problem & $m$ & $L$   & ACF &\%Setup&$\ttotal$\\ \hline \hline
(a) 5-point                  & $2095104$ & $19$ & $0.216$ & $29\%$ & $ 902$ \\ \hline
(b) 13-point $\kth{4}$ order & $4188160$ & $20$ & $0.262$ & $22\%$ & $1355$ \\ \hline
(c) Anis. rot. agnostic      & $4188162$ & $19$ & $\mathbf{0.816}$ & $ 5\%$ & $\mathbf{5453}$ \\ \hline
(d) Anis. rot. misaligned    & $3141633$ & $20$ & $\mathbf{0.870}$ & $ 4\%$ & $\mathbf{8136}$ \\ \hline
\end{tabular}
\caption{LAMG performance for grid graphs on a $1024 \times 1024$ grid with $n=1048576$ nodes.}
\label{results_negative}
\end{table}

LAMG exhibited mesh-independent convergence and run time in all cases and scaled linearly with grid size, albeit its convergence was much slower for cases (c) and (d), whose negative edge weights are more significant. Compared with the Bootstrap AMG method \cite{alg_distance_anis}, which focused on accurately finding the characteristic directions without sparing setup costs and only presented two-level experiments, LAMG is a full multi-level method with a far shorter setup time, although its ACF could also be significantly reduced using bootstrap tools. These results are certainly preliminary.

\subsection{Lean Geometric Multigrid}
\label{lmg}

Higher performance for the Poisson equation discretized on a uniform grid can be obtained by a standard 1:2 coarsening in every dimension at all levels and employing Gauss-Seidel relaxation in red-black ordering \cite[\S 3.6]{guide}. LAMG reduces to {\it Lean Geometric Multigrid}: standard multigrid cycle with index $\gamma=1.5$, first-order transfers and energy-corrected coarsening. Since the energy ratio is $2$ for all error modes, a flat correction $\mu=2$ is employed in (\ref{galerkin2}).

For the 2-D {\it periodic} Poisson problem, this cycle turns out to be a record-breaking Poisson solver in terms of asymptotic efficiency: it achieves a convergence factor of $.5$ per unit work, versus $.67$ for the classical multigrid V(1,1) cycle with linear interpolation and second-order full weighting \cite[\S 4.1]{lamg_arxiv}. For other boundary conditions, finding the right $\mu$ is not as easy; while supplementary local relaxations near boundaries theoretically ensure attaining the two-level rates \cite[\S 5]{guide}, it would be more beneficial to study the performance of adaptive energy correction in geometric LAMG.

\section{Future Research}
\label{extensions}
Enhancements and adaptations of the LAMG approach to related computational problems are outlined below.

\subsection{Coarsening Improvements}
\label{improvements}
The LAMG algorithm of \S \ref{algorithm} is by no means final and may be improved in various ways.
\bi
    \item The average setup time could be reduced by employing classical AMG with no test vectors, and switching to the LAMG strategy only when the former fails.
    \item In graphs with many weak edges (such as the outliers in Table.~\ref{outliers}), efficiency may be increased by temporarily ignoring them in the Galerkin operator computation, yet keeping track of their total contribution to each aggregate's stencil. If a level is reached at which this total is no longer small compared with the aggregate's other edge weights, it is reactivated.
    \item Currently, a node $u$ can only be aggregated with a direct neighbor $s$. In some problems, $u$'s second-degree neighbors should also be searched to ensure a good aggregation. For instance, in the anisotropic-rotated problem Fig.~\ref{negative_stencils}d, $u$ should be aggregated along the characteristic direction, i.e., with its southeast or northwest neighbor, neither of which is contained in $\cE_u$.
    \item If no small energy ratio can be found, or if subsequent cycle convergence is slow, isolated bottleneck nodes can be de-aggregated. Alternatively, one can up the interpolation caliber at these troublesome nodes, provided that this does not substantially increase the total coarse edges.
    \item Adaptive local relaxation sweeps may improve efficiency in various problems such as PDEs with structural singularities \cite{bai}.
\ei
Additionally, user-defined parameters could be supplied to treat special graph families more efficiently. For instance, if node coordinates are available, they can be used to generate smoother initial TVs than the default random initial guess. Optimizing the coarsening is most advantageous when $\bA \bx=\bb$ is solved for multiple $\bb$'s, since a larger setup cost is tolerable. Such is the case in time-dependent problems \cite{fischer}.

\subsection{Local Energy Correction}
\label{other_corrections}
Instead of a flat $\mu=\frac43$ factor in (\ref{galerkin2}), one can apply different $\mu$'s to different aggregates. We experimented with different energy correction schemes, some based on fitting the coarse nodal energies of TVs to their fine counterparts (cf.~(\ref{nodal_coarse})). While this can dramatically curtail energy inflation, care must be taken to avert over-fitting that ultimately results in the coarse-level correction operator's instability.

Analogously, one can define a local adaptive $\mu$ in the iterate recombination (\S \ref{adaptive}) at each level $l$, provided that it is properly smoothed \cite[\S 5.3]{lamg_arxiv}. It is unclear whether the extra work would be justified in either case. It may turn out useful for problems re-solved for many $\bb$ vectors, or when a larger setup overhead is tolerable.

\subsection{Other Linear Systems}
\label{other}

The LAMG caliber-1 aggregation can be applied to non-zero row sum matrices, except that the interpolation weights are no longer $1$. The affinity definition (\ref{cuv}) remains intact, while the corresponding $\bP$ entry is set to $p_{uv} := (X_u,X_v)/(X_v,X_v)$ (see (\ref{interp_accuracy})). Normally, relaxed TVs yield an accurate enough $p_{uv}$; in problems with almost-zero modes, e.g., the QCD gauge Laplacian, TVs may need to be improved by a bootstrap cycle \cite{bamg}.

Further research should be conducted for negative-weight graphs such as the high-order finite element and anisotropic grid graphs of \S \ref{negative_weights}. The reported convergence factors can be improved by producing bootstrapped TVs via applying multilevel cycles to $\bA \bx=\bzero$. The cycle is far more powerful than plain relaxation in damping smooth characteristic components, which should lead to more meaningful algebraic distances and to correct anisotropic coarsening in a second setup round (much larger spacing in the characteristic direction and no coarsening in the cross-characteristic direction). The bootstrap procedure should be useful in many other graphs.

\subsection{LAMG Eigensolver}
\label{eigenvalue}

The LAMG hierarchy can be combined with the Full Approximation Scheme (FAS) \cite[Chap.~8]{guide} to find the $K$ lowest eigenpairs of $\bA$, similarly to the work \cite{eis}. We perform the variable substitution $\bx^c = \bee^c + \bR \tbx$, transforming the coarse equation (\ref{galerkin}) into
\be
    \bA^c \bx^c = \bP^{T} \bb + \btau \tbx\,,\qquad \btau := \bA^c \bR - \bP^{T} \bA\,,
    \label{galerkin_fas}
\ee
followed by the fine-level correction $\tbx \leftarrow \tbx + \bP (\tbx^c - \bR \tbx)$. The elimination and aggregation are both special cases of (\ref{galerkin_fas}), with $\bR \tbx := \tbx_{\cC}$ and $\bR=\bzero$, respectively.

The analogue of (\ref{galerkin}) for coarsening $(\bA-\lambda_k \bI) \bx_k = \bzero$ is
\be
    \left( \bA^c - \lambda_k \bB^c \right) \bx^c_k = \btau_k \tbx\,,\qquad \btau_k \:=
    \left( \bA^c - \lambda_k \bB^c \right) \bR - \bP^{T} \left( \bA - \lambda_k \bB^c \right)\,.
    \label{galerkin_eigen}
\ee
Thus a separate affine term appears in the coarse equation of each approximate eigenvector $\bx^c_k$, $k=1,\dots,K$. In particular, the elimination of \S \ref{elimination} becomes approximate, yet (\ref{galerkin_eigen}) remains linear in $\lambda_k$, as opposed to the exact non-linear Schur complement formed by the AMLS method \cite{amls}. Gauss-Seidel may be replaced by Kaczmarz relaxation at very coarse levels to prevent the divergence of smooth error modes \cite{eis}.

Alternatively, one can incorporate the LAMG linear solver into a Rayleigh quotient iteration \cite[\S 8.2]{gvl},\cite{lobpcg}. However, FAS is attractive because it also applies to general nonlinear problems \cite[\S 8]{guide}, e.g., quadratic and linear programming.

\section{Conclusion}
Laplacian matrices underlie a plethora of graph computational applications ranging from genetic data clustering to social networks to fluid dynamics. To the best of our knowledge, the presented algorithm, Lean Algebraic Multigrid (LAMG), is the first graph Laplacian linear solver whose empirical performance approaches linear scaling for a wide variety of real-world graphs. Combinatorial Multigrid was also quite successful, performing faster on average, yet with many more outliers. The LAMG approach can also be generalized to non-diagonally-dominant, eigenvalue and nonlinear problems.

\section{Acknowledgments}
The authors wish to thank the referees for their fruitful comments, Ioannis Koutis, Tim Davis, David Gleich and Audrey Fu for useful discussions, Lorenzo Pesce for his help with porting LAMG to the Beagle Cray, and Dan Spielman for algorithmic discussions as well as \LaTeX\, typesetting advice.

\bibliographystyle{siam}       
\bibliography{lamg}

\end{document}

%% file: styles.tex
\usepackage{graphicx}
\usepackage{algorithm}
\usepackage{algorithmic}
\usepackage{amsmath}
\usepackage{amsfonts}
\usepackage{amssymb}
\usepackage{bm}
\usepackage{url}
\usepackage{multirow}

\newcommand{\ep}{\varepsilon}

\newcommand{\Etot}{E_{\mathrm{tot}}}

\newcommand{\numgraphs}{3774\;}
\newcommand{\maxedges}{47 million\;}
\newcommand{\setupmvm}{200}
\newcommand{\storageperedge}{4}
\newcommand{\solvemvm}{27}

\newcommand{\loget}{\log(1/\ep)}
\newcommand{\amax}{\alpha_{\mathrm{max}}}

\newcommand{\by}{\bm{\mathrm{y}}}
\newcommand{\bx}{\bm{\mathrm{x}}}

\newcommand{\bb}{\bm{\mathrm{b}}}
\newcommand{\bu}{\bm{\mathrm{1}}}

\newcommand{\bee}{\bm{\mathrm{e}}}
\newcommand{\br}{\bm{\mathrm{r}}}
\newcommand{\balpha}{\bm{\mathrm{\alpha}}}
\newcommand{\tbx}{\tilde{\bx}}

\newcommand{\tbee}{\tilde{\bee}}
\newcommand{\bzero}{\bm{\mathrm{0}}}

\newcommand{\bA}{\bm{\mathrm{A}}}
\newcommand{\bB}{\bm{\mathrm{B}}}

\newcommand{\bI}{\bm{\mathrm{I}}}
\newcommand{\bL}{\bm{\mathrm{L}}}

\newcommand{\bP}{\bm{\mathrm{P}}}

\newcommand{\bR}{\bm{\mathrm{R}}}
\newcommand{\bT}{\bm{\mathrm{T}}}

\newcommand{\btau}{\bm{\mathrm{\tau}}}
\newcommand{\bPi}{\bm{\mathrm{\Pi}}}

\newcommand{\cN}{\mathcal{N}}
\newcommand{\cE}{\mathcal{E}}

\newcommand{\cT}{\mathcal{T}}

\newcommand{\cF}{\mathcal{F}}
\newcommand{\cC}{\mathcal{C}}

\newcommand{\Real}{\mathbb{R}}

\newcommand{\matlab}{\textsc{Matlab\;}}

\newcommand{\agg}{\textsc{Aggregation\;}}

\newcommand{\tsetup}{t_{\mathrm{setup}}}
\newcommand{\tsolve}{t_{\mathrm{solve}}}
\newcommand{\ttotal}{t_{\mathrm{total}}}

\newcommand{\kth}[1]{#1^\mathrm{th}}
\newcommand{\lp}{\left(}
\newcommand{\rp}{\right)}

\newcommand{\mymin}[1]{\underset{#1}{\operatorname{min}}\;}
\newcommand{\mymax}[1]{\underset{#1}{\operatorname{max}}\;}
\newcommand{\argmin}[1]{\underset{#1}{\operatorname{argmin}}\;}

\newcommand{\be}{\begin{equation}}
\newcommand{\ee}{\end{equation}}
\newcommand{\bes}{\begin{equation*}}
\newcommand{\ees}{\end{equation*}}
\newcommand{\bea}{\begin{eqnarray}}
\newcommand{\eea}{\end{eqnarray}}
\newcommand{\beas}{\begin{eqnarray*}}
\newcommand{\eeas}{\end{eqnarray*}}
\newcommand{\bi}{\begin{itemize}}
\newcommand{\ei}{\end{itemize}}
\newcommand{\ben}{\begin{enumerate}}
\newcommand{\een}{\end{enumerate}}